\documentclass[12pt,twoside,a4paper]{article}
\title{\large{KELLER'S CONJECTURE ON THE EXISTENCE OF COLUMNS IN CUBE TILINGS OF $\er^n$}}
\date{}
\author{\\{Magdalena \L ysakowska and Krzysztof Przes{\l}awski}
\\
\small{Wydzia{\l} Matematyki, Informatyki i Ekonometrii, Uniwersytet Zielonog\'orski}\\
\small{ul. Z. Szafrana 4a, 65-516 Zielona G\'ora, Poland}\\
\small{M.Lysakowska@wmie.uz.zgora.pl}\\
\small{K.Przeslawski@wmie.uz.zgora.pl}}
\usepackage{amssymb,amsfonts,mathrsfs,amsmath}
\newtheorem{twr}{Theorem}

\newtheorem{lem}{Lemma}

\def\er{\mathbb{R}}
\def\zet{\mathbb{Z}}
\def\en{\mathbb{N}}

\begin{document}
\maketitle
\begin{abstract}
It is shown that if $n\le 6$, then  each tiling of $\er^n$ by translates of the unit cube $[0,1)^n$ contains a column; that is, a family  of the form $\{[0,1)^n+(s+ke_i)\colon k\in \zet\}$, where $s\in \er^n$ and $e_i$ is an element of the standard basis of $\er^n$. 

\end{abstract}
\section{Introduction}
\noindent
In his book \cite{Min}, which appeared in 1907, Hermann Minkowski proved
that every lattice tiling
of $\mathbb R^{n}$ by unit cubes contains two cubes that
have a common $(n-1)$-dimensional face,
whenever $n\leq 3$. This readily implies that
there is a column
of unit cubes contained in the tiling. On the other hand,
he conjectured
that the
same phenomenon holds in all dimensions.
In 1930, Otto Heinrich Keller \cite{K1}
extended
Minkowski's conjecture to
arbitrary
cube tilings of $\mathbb R^{n}$. In fact,
we have now two conjectures:
the
stronger, stated by Keller, which reads that each cube tiling
of $\mathbb R^{n}$ contains a column
and the weaker
which reads that each cube tiling of $\mathbb R^{n}$
contains
two
cubes that share an $(n-1)$-dimensional face.
In 1937, Keller published a short paper \cite{K2}
where he claimed that he proved his conjecture for $n\leq 6$. He also expressed a supposition that the conjecture
is not valid in dimensions
greater than 6. There are no  rigorous proofs in the paper
however.
In 1940, Oskar Perron \cite{P} published a complete proof of the weaker conjecture for dimensions not exceeding 6. 
He left aside the original 
conjecture of Keller. Apparently, he was focused
on
Minkowski's conjecture, which was verified by Haj\'os \cite{H} only one year later.
It was Perron who popularized the weaker conjecture under the name
of Keller.
In 1992, Jeff Lagarias and Peter Shor \cite{LS}
discovered a counterexample to the weaker conjecture
in dimension 10. Ten years later John Mackey \cite{Ma} found a
counterexample in dimension 8. This implies that both conjectures, the weaker
and the stronger,
are not valid in any dimension greater
than 7.  For dimension 7, both problems are 
completely open.
In the present note we show that the original Keller's conjecture is valid in all
dimensions up to 6. Our proof is based on Perron's
approach. However, certain substantial
modifications of
his method were necessary: He concentrated his attention on local
configurations of cubes, whereas we have to play with all cubes of tiling.
This 
forces us to work with 
an appropriate enumeration
of cubes.

\section{The existence of columns}
\noindent
We define a \textit{cube}
in the $n$-dimensional Euclidean space $\er^n$ to
be any translate of the unit cube $[0,1)^n$. Let $T$ be a subset of $\er^n$. The family
$[0,1)^{n}+T:= \{[0,1)^{n}+t\colon t\in T\}$ is said to be  a \textit{cube tiling}
of $\mathbb R^{n}$ if for each
pair of distinct vectors $s,t\in T$ the cubes $[0,1)^{n}+s$
and $[0,1)^{n}+t$ are disjoint and
$\bigcup[0,1)^{n}+T=\mathbb R^{n}$. We refer to  $T$
as a set that determines a cube tiling.

As usual, we denote by $\zet$ the set of all
integers while
the set of positive integers is denoted by $\en$. Let $n\in \en$. The set $\{1,2,\ldots,n\}$ is denoted by $[n]$.  

Let us recall the following fundamental result from \cite{K1}.
\begin{twr}[O. H. Keller, 1930]
If $[0,1)^{n}+T$ is a cube tiling of $\mathbb{R}^{n}$,
then for each pair of distinct elements $s,t\in T$ there is $j\in[n]$ such that
$|s_{j}-t_{j}|\in\mathbb{N}$.
\end{twr}

\noindent
{\bf Proof}.
Suppose the theorem
is not valid. Then there is a set $T$ that determines  a
cube tiling of $\mathbb R^{n}$
which  contains a pair of distinct
elements $s,t\in T$ such that 
$s_{j}-t_{j}
\notin\mathbb{Z}\setminus\{0\}$ for each $j\in[n]$. Clearly, the set $S=T-t$ determines a cube
tiling as well.
Let $u=s-t$. The elements $u$ and $0$ belong to $S$ and
$u_j\notin\mathbb{Z}\setminus\{0\}$ for each $j\in [n]$.
For $x\in \er^n$, let
$i(x):=\operatorname{max}\{i\colon |x_{j}|<1\,\, \text{for}\,  j\leq i\}$, where
we have assumed that $\operatorname{\max}\emptyset=0$. We have $i(u)<n$, otherwise
cubes $[0,1)^{n}$ and $[0,1)^{n}+u$
would intersect contradicting the assumption that
$S$
determines a  cube tiling.
Let $k:=i(u)+1$ and let
$V:=\{v\in S\colon v_{k}-u_{k}\in\mathbb{Z}\}$. If $\ell$ is a straight line intersecting
one of the cubes $[0,1)^{n}+v$, $v\in V$, which in addition is parallel to the $k$-th
coordinate axis, then
$l\subset\bigcup_{v\in V}{([0,1)^{n}+v)}$. This observation
leads to the conclusion that the set $U:=(V-\lfloor u_{k}\rfloor e_{k})\cup(S\setminus V)$
determines a cube tiling. Obviously, $r:=u-\lfloor u_{k}\rfloor e_{k}$
and $0$ are elements
of $U$.
Moreover,  $i(r)=i(u)+1$ and  $r_{j}\notin\mathbb{Z}\setminus\{0\}$
for every $j\in[n]$.
Now, we can replace $S$ by $U$, and continue in this manner
eventually arriving to a set that determines a cube tiling and contains $0$ and an
element $w$ such
that $i(w)=n$ and $w_{j}\notin\mathbb{Z}\setminus\{0\}$ for
every $j\in[n]$ which, as we know, is impossible.\hfill $\square$\\

Let us suppose that for each $j\in[n]$ a mapping $\varepsilon_{j}\colon\mathbb R\to\mathbb N$
is given
such that
for each element $x\in\mathbb R$ the restriction
$\varepsilon_{j}|x+\mathbb Z$ is a bijection between the sets
$x+\mathbb Z:=\{x+k\colon k\in\mathbb Z\}$ and $\mathbb N$. The mapping
$\varepsilon\colon\mathbb R^{n}\to\mathbb N^{n}$ defined by the formula
$$
\varepsilon (x)=\varepsilon(x_{1},\ldots,x_{n})=(\varepsilon_{1}(x_{1}),\ldots,\varepsilon_{n}(x_{n}))
$$
is said to be a \textit{natural code} (of $\er^n$). The vector $\varepsilon(x)$
is referred to as the code of $x$.

\begin{twr}\label{kod}
Fix a natural
code
$\varepsilon\colon\mathbb R^{n}\to\mathbb N^{n}$. Then a set $T\subseteq\mathbb R^{n}$
determines a cube tiling of $\mathbb R^{n}$ if and only if  $\varepsilon(T)=\en^n$ and
for every pair of distinct elements $s,t\in T$ there is $j\in [n]$ such
that
$|s_j-t_j|\in \en$.
\end{twr}

\noindent
{\bf Proof}.
\noindent $(\Rightarrow)$
By Keller's theorem, it suffices to show that
$\varepsilon(T)=\en^n$. We proceed by induction with respect $n$. For $n=1$ the 
assertion is a consequence of the definition of
a natural code and the fact
that each set determining a cube tiling of $\er^1$  coincides with one of the cosets
$\zet +x$, $x\in \er$. Fix $k\in \en$, and define the set
$$
T^{k}:=\{t=(t_{1},\ldots,t_{n})\in T\colon \varepsilon_{n}(t_{n})=k\}.
$$
If
$n>1$, then $\er^{n-1}$ is a non-trivial Euclidean space. Let $T_{n'}^k\subset \er^{n-1}$
be the image of $T^k$ under the projection
$x=(x_1,\ldots,x_n)\mapsto x_{n'}:=(x_1,\ldots,x_{n-1})$ of $\er^{n}$ onto $\er^{n-1}$.
Let us show now that $T_{n'}^k$
determines a cube tiling of
$\er^{n-1}$. Let $x\in\mathbb R^{n-1}$ and let $T(x)=\{t\in T\colon ([0,1)^{n}+t)\cap(\{x\}\times\mathbb R)\neq\emptyset\}$.
Since $[0,1)^{n}+T$ is a cube tiling of $\mathbb R^{n}$,   the set $T(x)+[0,1)^n$
is a covering of
the set $\{x\}\times\mathbb R$ by disjoint sets. Therefore,
$T(x)_{n}=\{t_{n}\colon t\in T(x)\}$ determines a cube tiling of $\mathbb R$.
The latter set is transformed  by $\varepsilon_{n}$ bijectively onto $\mathbb N$. In
particular, there is an element
$t\in T(x)$ such that $\varepsilon_{n}(t_{n})=k$;
equivalently,  $t\in T^{k}$ and $x\in [0,1)^{n-1}+t_{n'}$.
Since $x$ is arbitrary,
we deduce that $T^k_{n'}$ determines a cube tiling of $\er^{n-1}$, as announced.
Let
$\varepsilon_{n'}:=(\varepsilon_1,\ldots, \varepsilon_{n-1})$. By the induction
hypothesis,
$\varepsilon_{n'}(T^k_{n'})=\en^{n-1}$. Therefore,
$\varepsilon(T^k)=\en^{n-1}\times\{k\}$. Finally,
$$
\varepsilon(T)=\bigcup_{k\in \en}\varepsilon(T^k)=\bigcup_{k\in \en} \en^{n-1}\times\{k\}=\en^n.
$$
\noindent$(\Leftarrow)$
As $\varepsilon$ maps $T$
`onto' $\mathbb N^{n}$, it maps $T^{k}$ `onto'
$\mathbb N^{n-1}\times\{k\}$.
Hence $\varepsilon_{n'}$
maps  $T_{n'}^{k}$ `onto'
$\mathbb N^{n-1}$. Moreover, by our assumptions
for every two elements $s,t\in T^{k}_{n'}$
there is $j\in[n-1]$ such that
$|s_{j}-t_{j}|\in\mathbb{N}$. Therefore,  $T_{n'}^{k}$
determines a cube tiling of $\er^{n-1}$ by the induction
 hypothesis.
Let $x\in\mathbb R^{n-1}$.
Then for each number $k\in\mathbb N$ there is $t^{k}\in T^{k}$ such that
$x\in[0,1)^{n-1}+ t^k_{n'}$. Since for every pair of distinct elements
$r,s\in \{t^{k}\colon k\in\mathbb N\}$ there is $j\in[n]$ such that
$|r_{j}-s_{j}|\in\mathbb{N}$
and the cubes $[0,1)^{n-1}+ t^k_{n'}$, $k\in \en$, intersect,
the set
$\{t^{k}_{n}\colon k\in\mathbb N\}$ determines a cube tiling of $\mathbb R$.
Thus,
$\{x\}\times\mathbb R\subseteq\bigcup_{k\in\mathbb N}([0,1)^{n}+t^{k})\subseteq\bigcup([0,1)^{n}+T)$.
Consequently,
$\bigcup([0,1)^{n}+T)=\mathbb R^{n}$. \hfill$\square$\\

Two vectors $x$ and $y\in\mathbb R^{n}$ are 
\textit{$\zet$-distinguishable}
if there is $j\in[n]$ such that $x_{j}\neq y_{j}$
and $y_{j}\in x_{j}+\mathbb Z$. A system of vectors is called $\zet$-distinguishable, or shortly distinguishable, if any two vectors  of this system are $\zet$-distinguishable.

Let $l<n$. A family of boxes $F$ is said to be an \textit{$l$-column}
if there is a set of vectors $S\subset\mathbb R^{n}$ such that the
following conditions
are satisfied:
\begin{itemize}
\item[(1)] $F=S+[0,1)^n$;
\item[(2)] there is $i\in [n]$
such that the mapping $x\mapsto x_i$ transforms $S$ bijectively into a set that determines
a cube partition of $\er^1$;
\item[(3)]there are $l$
indices $j\in [n]$ such that the
sets $S_j:= \{x_j: x\in S\}$ are singletons.
\end{itemize}
We refer to $S$ as a set that determines an $l$-column.
Every $(n-1)$-column contained in $\mathbb R^{n}$ is called 
in short a \textit{ column}.

Two sets of vectors $F,G\subseteq\mathbb R^{n}$  are
{\it isomorphic}, if there are a
bijection $f\colon F\to G$ and a permutation $\sigma\colon[n]\to[n]$ such that
for every
pair of vectors $x,y\in F$ and for each
 $j\in[n]$
the following conditions are satisfied:
\begin{itemize}
\item[(1)] $x_{j}=y_{j}$
if and only if $f(x)_{\sigma(j)}=f(y)_{\sigma(j)}$,
\item[(2)]  $|x_{j}-y_{j}|\in\mathbb N$ if and only if
 $|f(x)_{\sigma(j)}-f(y)_{\sigma(j)}|\in\mathbb N$.
\end{itemize}

Let $x$ be an element of $\er^n$. In what follows, we often write
$x: x_{1}\ldots x_{n}$ instead of  $x=(x_{1},\ldots,x_{n})$,  according to the convention adopted by Perron \cite{P}.
The coordinates of the vectors belonging to any set that determines a cube tiling of $\er^n$ are denoted by Roman or Greek lower case
letters. We apply somewhat untypical convention that if various coordinates of a vector are represented by the same letter (possibly with the same upper and lower indices), then it does not imply that they have the same value. Starting from the proof of Theorem \ref{ge4}, it is tacitly assumed that
when we talk about the set of vectors that determines a cube tiling, a natural code of $\er^n$ is already defined. Such a code serves as a system of coordinates of the cube tiling. Lower indices of coordinates of a vector correspond to the code of this
vector; e.g.  the string of symbols $w: a_{1} \alpha_{3}^5 a_{2}$ means that the vector $w=(a_{1}, \alpha_{3}^5, a_{2})$ has the code $\varepsilon(w)=(1,3,2)$. 
If we have two vectors 
whose $i$-th coordinates are denoted by the same lower case letter
 with the same upper index, if there is any, then they differ by an integer; e.g. if $w_{1}: a_{1} a_{1} \alpha^2_{1}$
and $w_{2}: a_{2} b_{1} \alpha^2_{3}$, then $(w_1)_1-(w_2)_1=a_1-a_2$ and $(w_1)_3-(w_2)_3=\alpha^2_1-\alpha^2_3$
are non-zero integers.
 If $i$-th coordinates of two vectors are denoted by different Roman letters, then their difference is not an integer; e.g. if $w_1$ and $w_2$ are as above, then $(w_1)_2=a_1$ while $(w_2)_2=b_1$, therefore, 
$(w_1)_2-(w_2)_2=a_1-b_1$ is not an integer.
As is seen from the examples, 
Roman letters occur only with lower indices while Greek letters are equipped with lower and upper indices. The value of a coordinate of a vector is denoted by a Greek letter when it is not explicit; e.g.  the third
coordinates of the vectors $w_{1}: a_{1} a_{2} \alpha_{2}^{3}$ and $w_{2}: a_{1} a_{1} \alpha_{2}^{4}$ have the same code but it is not decided 
whether they are equal or different. There is only one place, the proof of Lemma \ref{le6}, where coordinates of a
vector are denoted by Roman capital letters. Their use will be explained therein. 

\begin{twr}\label{sys}
Let $n$ be a positive integer. If every  cube tiling of $\mathbb R^{n}$
contains a column, then every  cube tiling of $\mathbb R^{m}$ contains an
$(n-1)$-column, for each $m>n$. 
\end{twr}

\begin{twr}\label{ge2}
Every  cube tiling of $\mathbb R^{2}$ contains a column.
\end{twr}

\noindent
{\bf Proof}.
Let $T$ be an arbitrary set which determines a cube tiling of $\mathbb R^{2}$ and let
$\varepsilon\colon\mathbb R^{2}\to\mathbb N^{2}$ be a natural code.
Let us consider these  elements of $T$ which have the codes $(k,1)$, $k\in\mathbb N$. By Theorem \ref{kod} and the notation introduced above, they can be written as follows
$$
\begin{array}{lllr}
w_{1,\;1}: & a_{1} & a_{1}, \\
\\
w_{1,\;l}: & a_{l} & \alpha_{1}^{l-1},& \qquad l\geq2.
\end{array}
$$
If  $\alpha_{1}^{i}=a_{1}$ for all $i\geq1$, then the vectors $w_{1,k},\ k\geq1$, determine a column.
Let us suppose now that at least one of the numbers $\alpha_{1}^{i},\ i\geq1$, is distinct from $a_{1}$.
We can assume that
it is $\alpha_{1}^{1}$, as if $\alpha_{1}^{1}=a_{1}$ and for example $\alpha_{1}^{5}\neq a_{1}$,
then we can change our natural
code replacing $\varepsilon_1$ by the composite $\tau\circ\varepsilon_1$, where $\tau$ is the transposition $(2\,6)$, arriving in this way to a system of vectors where $w_{1,\;2}$ has its second coordinate different from $a_1$. Now, let us take into account
all  vectors belonging to $T$ with codes $(1,l)$, $l\geq2$.
As  these vectors together with $w_{1,\;1}$ and $w_{1,\;2}$ are $\zet$-distinguishable, they must have the following form
$$
\begin{array}{lllr}
w_{2,\;l-1}: & a_{1} & a_{l}, & \qquad l\geq2.\\
\end{array}
$$
Thus, they together with $w_{1,\;1}$ determine a column.
\hfill$\square$

\begin{twr}\label{ge3}
Every cube tiling of $\mathbb R^{3}$ contains a column.
\end{twr}

\noindent
{\bf Proof}.
Let $T$ be an arbitrary set which determines a cube tiling of $\mathbb R^{3}$ and let $\varepsilon\colon\mathbb R^{3}\to\mathbb N^{3}$ be a natural code.
By Theorems \ref{sys} and \ref{ge2}, the set $T$ contains vectors which
determine a $1$-column.
Passing
to an isomorphic system if necessary, we can assume that the vectors determining our $1$-column are as follows
$$
\begin{array}{llllr}
w_{1,\;1}: & a_{1} & a_{1} & a_{1},\\
\\
w_{1,\;l}: & a_{l} & a_{1} & \alpha_{1}^{l-1}, & \qquad l\geq2.
\end{array}
$$
If  $\alpha_{1}^{i}=a_{1}$ for all $i\geq1$, then the vectors $w_{1,\;k}$, $k\geq1$, determine a column.
Suppose that $\alpha_{1}^{i}\neq a_{1}$ for some $i\geq1$.
As previously, we can assume that
$\alpha_{1}^{1}\neq a_{1}$. Consider all vectors belonging to $T$ with codes $(1,1,l)$, $l\geq2$.
As these vectors together with $w_{1,\;1}$ and $w_{1,\;2}$ are $\zet$-distinguishable, they can be written in the following
form
$$
\begin{array}{llllr}
w_{2,\;l-1}: & a_{1} & \beta_{1}^{l-1} & a_{l}, & \qquad l\geq2.
\end{array}
$$
If  $\beta_{1}^{i}=a_{1}$ for all $i\geq1$, then the vectors $w_{1,\;1}$, $w_{2,\;k}$, $k\geq1$, determine a column.
Suppose that at least one of the coordinates $\beta_{1}^{i}$, $i\geq1$, is different from $a_{1}$. By the same reason as in the preceding proof, we can assume that
$\beta_{1}^{1}\neq a_{1}$. Now, let us consider all vectors from $T$ with codes
$(1,l,2)$, $l\geq2$.
As  these vectors together with $w_{1,\;1}$, $w_{1,\;2}$ and $w_{2,\;1}$ are $\zet$-distinguishable, they must have the following form
$$
\begin{array}{llllr}
w_{3,\;l-1}: & a_{1} & \beta_{l}^{1} & a_{2}, & \qquad l \geq 2.
\end{array}
$$
Thus, they and $w_{2,\;1}$ determine a column.
\hfill$\square$

\begin{twr}\label{ge4}
Every cube tiling of $\mathbb R^{4}$ contains a column.
\end{twr}

\noindent
{\bf Proof}.
Let $T$ be an arbitrary set which determines a cube tiling of $\mathbb R^{4}$.
By Theorems \ref{sys} and \ref{ge3}, the set $T$
contains vectors which
determine a $2$-column.
Passing
to an isomorphic system if  necessary, we can assume that the vectors determining our $2$-column are as follows
$$
\begin{array}{lllllr}
w_{1,\;1}: & a_{1} & a_{1} & a_{1} & a_{1},\\
\\
w_{1,\;l}: & a_{l} & a_{1} & a_{1} & \alpha_{1}^{l-1},& \qquad l\geq2.
\end{array}
$$
If  $\alpha_{1}^{i}=a_{1}$ for all $i\geq1$, then the vectors $w_{1,\;k}$, $k\geq1$, determine a column.
Suppose that $\alpha_{1}^{i}\neq a_{1}$ for some $i\geq1$. We can assume that
$\alpha_{1}^{1}\neq a_{1}$.  Now, let us consider all vectors from $T$ with codes
$(1,\ast,\ast,l)$, $l\geq2$, where $\ast$ can take any value
from $\mathbb N$, such that
their second and third coordinates
are different from $a_{l}$, $l\geq2$. (Such vectors exist, e.g. 
the vector whose code is $(1,1,1,5)$  has the required property.)
Let us pick a vector from among them whose \textit{middle} coordinates (second and third) differ from $a_1$ at
as many places as possible. 
Let us change the natural code so that the vectors $w_{1,\;l}$, $l\ge 1$ remain unaffected while the picked vector has its code equal to $(1,1,1,2)$. By $\zet$-distinguishability of this vector from $w_{1,\;1}$ and $w_{1,\;2}$, it can be written as follows
$$
\begin{array}{lllll}
w_{2,\;1}: & a_{1} & \beta_{1}^{1} & \beta_{1}^{2} & a_{2}.
\end{array}
$$ 
Three cases have to be considered:

\medskip
\noindent
\textit{Case 1.} $\beta_{1}^{1}=\beta_{1}^{2}=a_{1}$.

\medskip
Then take all  the vectors with codes  $(1,1,1,l)$, $l \ge 3$.  Since each of them must be distinguishable from the vectors  
$w_{1,\;1}$, $w_{1,\;2}$ and $w_{2,\;1}$, we
deduce that they can be written in the form 

$$
\begin{array}{lllllr}
w''_{2,\;l-1}: & a_{1} & \beta_{1}^{2l-3} & \beta_{1}^{2l-2} & a_{l}, & \qquad l\geq3.
\end{array}
$$
As $w_{2,\;1}$ has the smallest possible number of
the middle coordinates equal to $a_{1}$, we have   
$\beta_{1}^{i}=a_{1}$ for all $i\geq3$. Thus, the vectors
$w_{1,\;1}$, $w_{2,\;1}$, $w''_{2,\;l}$, $l\geq2$, determine a column.

\medskip
\noindent
\textit{Case 2.} Exactly one of the coordinates
$\beta_{1}^{1}$ and $\beta_{1}^{2}$ is equal
to $a_{1}$.  

\medskip
Then we can assume that $\beta_{1}^{2}=a_{1}$, as in the other case we would change the order of the second and the third coordinate. Consider all vectors with codes
$(1,l,1,2)$,  $l \ge 2$. By
the distinguishability,  they can be written in the form
$$
\begin{array}{lllllr}
w'_{2,\;l}: & a_{1} & \beta_{l}^{1} & \gamma_{1}^{l-1} & a_{2}, & \qquad l\geq2.
\end{array}
$$
As $w_{2,\;1}$ has the smallest possible number of
the middle coordinates equal to $a_{1}$, we have
$\gamma_{1}^{i}=a_{1}$ for all $i\geq1$. Thus, the vectors
$w_{2,\;1}$ and $w'_{2,\;l}$, $l\geq2$, determine a column.

\medskip
\noindent
\textit{Case 3.}
Both coordinates $\beta_{1}^{1}$ and $\beta_{1}^{2}$ are different from $a_{1}$.
\medskip

Then consider all  vectors with codes $(1,l,1,2)$ and $(1,1,l,2)$, $l\geq2.$
By their distinguishability from $w_{1,\;1}$, $w_{1,\;2}$ and $w_{2,\;1}$, they have the form
$$
\begin{array}{lllllr}
w_{2,\;l}: & a_{1} & \beta_{l}^{1} & \gamma_{1}^{l-1} & a_{2},\\
\\
w_{3,\;l-1}: & a_{1} & \delta_{1}^{l-1} & \beta_{l}^{2} & a_{2}, & \qquad l\geq2. 
\end{array}
$$
Since the vectors $w_{2,\;l}$, $l\geq2$, and $w_{3,\;k}$, $k\geq1$, are $\zet$-distinguishable, we have $\delta_{1}^{i}=\beta_{1}^{1}$ for all $i\geq1$
or $\gamma_{1}^{i}=\beta_{1}^{2}$ for all $i\geq1$. If the first case takes place the vectors $w_{2,\;1}$, $w_{3,\;k}$, $k\geq1$, determine a column; otherwise the vectors $w_{2,\;l}$, $l\geq1$, determine a column.
\hfill$\square$

\begin{twr}\label{ge5}
Every cube
tiling of $\mathbb R^{5}$ contains a column.
\end{twr}

\noindent
{\bf Proof}.
Let $T$ be an arbitrary set which determines a cube tiling of $\mathbb R^{5}$.
By Theorems \ref{sys} and \ref{ge4}, the set $T$ contains vectors which
determine a $3$-column.
Passing
to an isomorphic system if  necessary, we can suppose that the vectors determining our $3$-column are as follows
$$
\begin{array}{llllllr}
w_{1,\;1}: & a_{1} & a_{1} & a_{1} & a_{1} & a_{1},\\
\\
w_{1,\;l}: & a_{l} & a_{1} & a_{1} & a_{1} & \alpha_{1}^{l-1}, & \qquad l\geq2.
\end{array}
$$
If  $\alpha_{1}^{i}=a_{1}$ for all $i\geq1$, then
the vectors $w_{1,\;k}$, $k\geq1$, determine a column. Suppose that $\alpha_{1}^{i}\neq a_{1}$ for some $i\geq1$. As before, we can assume that
$\alpha_{1}^{1}\neq a_{1}$. Now, let us consider all  vectors from $T$ with codes
$(1,\ast,\ast,\ast,l)$, $l\geq2$,
where $\ast$ can take  any value
from $\mathbb N$, such  that their middle coordinates (second, third and fourth)
are different from $a_{l}$, $l\geq2$. Let us pick a vector from among them  whose middle coordinates differ from  $a_{1}$ at as many places as possible.
Similarly as in the proof of Theorem \ref{ge4}, we can assume that this vector has its code equal to $(1,1,1,1,2)$. By its distinguishability from $w_{1,\;1}$ and $w_{1,\;2}$, it can be written in the form
$$
\begin{array}{llllll}
w_{2,\;1}: & a_{1} & \beta_{1}^{1} & \beta_{1}^{2} & \beta_{1}^{3} & a_{2}. 
\end{array}
$$

\noindent
Four cases have to be considered:

\medskip
\noindent
\textit{Case 1.}
$\beta_{1}^{1}=\beta_{1}^{2}=\beta_{1}^{3}=a_{1}$.

\medskip
Then take all vectors with codes
$(1,1,1,1,l)$, $l \geq 3$. By the distingui\-sha\-bility, these vectors have the form
$$
\begin{array}{llllllr}
w'''_{2,\;l-1}: & a_{1} & \beta_{1}^{3l-5} & \beta_{1}^{3l-4} & \beta_{1}^{3l-3} & a_{l},  & \qquad l\geq3.\\
\end{array}
$$
As $w_{2,\;1}$ has the smallest possible number  of the middle coordinates equal to $a_{1}$, we have $\beta_{1}^{i}=a_{1}$ for all $i\geq4$. Thus, the vectors $w_{1,\;1}$, $w_{2,\;1}$, $w'''_{2,\;l}$, $l\geq2$, determine a column.

\medskip
\noindent
\textit{Case 2.}
Exactly one of the coordinates  $\beta_{1}^{1}$, $\beta_{1}^{2}$, $\beta_{1}^{3}$
is distinct from $a_{1}$.
\medskip

Then we can assume that $\beta_{1}^{1}\neq a_{1}$, $\beta_{1}^{2}=\beta_{1}^{3}=a_{1}$, as in the other case
we would change the order of the appropriate coordinates.  Take into account all vectors  with codes $(1,l,1,1,2)$, $l\geq2$. They have the following form
$$
\begin{array}{llllllr}
w''_{2,\;l}: & a_{1} & \beta_{l}^{1} & \gamma_{1}^{2l-3} & \gamma_{1}^{2l-2} & a_{2}, & \qquad l\geq2.
\end{array}
$$
As $w_{2,\;1}$ has the smallest possible number of the middle coordinates equal to $a_{1}$, we have $\gamma_{1}^{i}=a_{1}$ for all $i\geq1$.
Thus, the vectors $w_{2,\;1}$, $w''_{2,\;l}$, $l\geq2$, determine a column.

\medskip
\noindent
\textit{Case 3.}
Exactly two of the coordinates  $\beta_{1}^{1}$, $\beta_{1}^{2}$, $\beta_{1}^{3}$
are distinct from $a_{1}$.
\medskip

Then, by the same reason as before, we can assume that $\beta_{1}^{1}\neq a_{1}$, $\beta_{1}^{2}\neq a_{1}$, and $\beta_{1}^{3}=a_{1}$. Take into account all  vectors with codes $(1,l,1,1,2)$ and $(1,1,l,1,2)$, $l\geq2$.
As these vectors together with $w_{1,\;1}$, $w_{1,\;2}$ and $w_{2,\;1}$  are distinguishable, they can be written in the following form
$$
\begin{array}{llllllr}
w'_{2,\;l}: & a_{1} & \beta_{l}^{1} & \gamma_{1}^{2l-3} & \gamma_{1}^{2l-2} & a_{2},\\
\\
w'_{3,\;l-1}: & a_{1} & \delta_{1}^{2l-3} & \beta_{l}^{2} & \delta_{1}^{2l-2} & a_{2}, & \qquad l\geq2. 
\end{array}
$$
Since the vectors $w'_{2,\;l}$, $l\geq2$, and $w'_{3,\;k}$, $k\geq1$,
must be distinguishable, we have $\delta_{1}^{i}=\beta_{1}^{1}$ for  $i=1,3,5,\ldots$
or $\gamma_{1}^{i}=\beta_{1}^{2}$ for  $i=1,3,5,\ldots$. If  the first possibility happens,
then, as $w_{2,\;1}$ has the smallest possible number of the middle coordinates equal to $a_{1}$, we have $\delta_{1}^{i}=a_{1}$ for  $i=2,4,\ldots$. Thus, the vectors $w_{2,\;1}$, $w'_{3,\;k}$, $k\geq1$, determine a column. If the  second possibility takes place, then we have
$\gamma_{1}^{i}=a_{1}$ for  $i=2,4,\ldots$, and the
vectors $w_{2,\;1}$, $w'_{2,\;k}$, $k\geq2$,
determine a column.

\medskip
\noindent
\textit{Case 4.}
All of the coordinates  $\beta_{1}^{1}$, $\beta_{1}^{2}$, $\beta_{1}^{3}$
are distinct from $a_{1}$.
\medskip

Then consider all vectors belonging to $T$ with codes $(1,l,1,1,2)$, $l\geq2$.
Since each of these vectors is distinguishable from $w_{1,\;1}$, $w_{1,\;2}$ and $w_{2,\;1}$, they have the form
$$
\begin{array}{llllllr}
w_{2,\;l}: & a_{1} & \beta_{l}^{1} & \gamma_{1}^{2l-3} & \gamma_{1}^{2l-2} & a_{2}, & \qquad l\geq2. 
\end{array}
$$
If $\gamma_{1}^{i}=\beta_{1}^{2}$ for  $i=1,3,\ldots$ and
$\gamma_{1}^{i}=\beta_{1}^{3}$
for  $i=2,4,\ldots$, then the vectors $w_{2,\;l}$, $l\geq1$,
determine a column. Suppose that at least one of the above equalities  does not happen.
We can assume that $\gamma_{1}^{1}\neq\beta_{1}^{2}$, as if $\gamma_{1}^{1}=\beta_{1}^{2}$ and for example
$\gamma_{1}^{4}\neq\beta_{1}^{3}$,
then we would change the order of the third and forth coordinates, and replace the code $\varepsilon$ with $\varepsilon':=(\varepsilon_1, \tau\circ\varepsilon_2, \varepsilon_4,\varepsilon_3,\varepsilon_5)$, where $\tau$ is the transposition $(2\,3)$.  
Take all vectors  with codes
$(1,1,l,1,2)$, $l\geq2$. They can be written in the form
$$
\begin{array}{llllllr}
w_{3,\;l-1}: & a_{1} & \delta_{1}^{2l-3} & \beta_{l}^{2} & \delta_{1}^{2l-2} & a_{2},  & \qquad l\geq2. 
\end{array}
$$
The distinguishability of  the vectors $w_{2,\;2}$ and
$w_{3,\;k}$, $k\geq1$ implies 
$\delta_{1}^{i}=\beta_{1}^{1}$ for  $i=1,3,\ldots$, as $\gamma_{1}^{1}\neq\beta_{1}^{2}$.
If now  $\delta_{1}^{i}=\beta_{1}^{3}$ for  $i=2,4,\ldots$, then the vectors
$w_{2,\;1}$, $w_{3,\;l}$, $l\geq1$, determine a column. Therefore, we can
assume that at least one of the coordinates
$\delta_{1}^{i}$, $i=2,4,\ldots$, is distinct from $\beta_{1}^{3}$.
We can also assume that $\delta_{1}^{2}\neq\beta_{1}^{3}$, as if  $\delta_{1}^{2}=\beta_{1}^{3}$
and for example $\delta_{1}^{4}\neq\beta_{1}^{3}$,
then we would change the code replacing $\varepsilon_3$ by the composite $\tau\circ \varepsilon_3$, where $\tau=(2\,3)$. 
Now, let us take into account all  vectors with codes $(1,1,1,l,2)$, $l\geq2$.
Since each of them must be distinguishable from $w_{1,\;1}$, $w_{1,\;2}$ and
$w_{2,\;1}$, they can be written as follows
$$
\begin{array}{llllllr}
w_{4,\;l-1}: & a_{1} & \varepsilon_{1}^{2l-3} & \varepsilon_{1}^{2l-2} & \beta_{l}^{3} & a_{2}, &\qquad l\geq2. 
\end{array}
$$
Since $\delta_{1}^{2}\neq\beta_{1}^{3}$,
by the distinguishability of the vectors $w_{3,\;1}$ and $w_{4,\;l}$, $l\geq1$,
we have $\varepsilon_{1}^{i}=\beta_{1}^{2}$ for  $i=2,4,\ldots$. Again, as the vectors
$w_{4,\;l}$, $l\geq1$, and  $w_{2,\;2}$ are
distinguishable, $\varepsilon_{1}^{i}=\beta_{1}^{1}$ for
$i=1,3,\ldots$
or $\gamma_{1}^{2}=\beta_{1}^{3}$. If the first possibility takes place, then the
vectors $w_{2,\;1}$, $w_{4,\;l}$, $l\geq1$, determine a column. Hence, the second possibility  has to be considered. Consequently, we obtain
$$
\begin{array}{llllllll}
w_{1,\;1}: & a_{1} & a_{1} & a_{1} & a_{1} & a_{1},\\
\\
w_{1,\;2}: & a_{2} & a_{1} & a_{1} & a_{1} & \alpha_{1}^{1},\\
\\
w_{2,\;1}: & a_{1} & \beta_{1}^{1} & \beta_{1}^{2} & \beta_{1}^{3} & a_{2},\\
\\
w_{2,\;2}: & a_{1} & \beta_{2}^{1} & \gamma_{1}^{1} & \beta_{1}^{3} & a_{2},\\
\\
w_{3,\;1}: & a_{1} & \beta_{1}^{1} & \beta_{2}^{2} & \delta_{1}^{2} & a_{2},\\
\\
w_{4,\;1}: & a_{1} & \varepsilon_{1}^{1} & \beta_{1}^{2} & \beta_{2}^{3} & a_{2}.
\end{array}
$$
Take into account all vectors with codes $(1,1,2,l,2)$, $(1,l,1,2,2)$, $(1,2,l,1,2)$, $l\geq2$.
They have the form
$$
\begin{array}{llllllr}
w_{5,\;l-1}: & \eta_{1}^{2l-3} & \beta_{1}^{1} & \beta_{2}^{2} & \delta_{l}^{2} & \eta_{2}^{2l-2},\\
\\
w_{6,\;l-1}: & \mu_{1}^{2l-3} & \varepsilon_{l}^{1} & \beta_{1}^{2} & \beta_{2}^{3} & \mu_{2}^{2l-2},\\
\\
w_{7,\;l-1}: & \nu_{1}^{2l-3} & \beta_{2}^{1} & \gamma_{l}^{1} & \beta_{1}^{3} & \nu_{2}^{2l-2}, & \qquad l\geq2. 
\end{array}
$$
If $\eta_{1}^{i}=a_{1}$ for  $i=1,3\ldots$ and $\eta_{2}^{i}=a_{2}$ for  $i=2,4\ldots$, then the
vectors
$w_{3,\;1}$, $w_{5,\;l}$, $l\geq1$, determine a column.
If $\mu_{1}^{i}=a_{1}$ for  $i=1,3\ldots$ and $\mu_{2}^{i}=a_{2}$ for  $i=2,4\ldots$, then the vectors
$w_{4,\;1}$, $w_{6,\;l}$, $l\geq1$, determine a column.
If $\nu_{1}^{i}=a_{1}$ for  $i=1,3\ldots$ and
$\nu_{2}^{i}=a_{2}$ for  $i=2,4\ldots$, then the vectors
$w_{2,\;2}$, $w_{7,\;l}$, $l\geq1$, determine a column.
Therefore, we can  assume that each of the above three statements is false. Then, by the distinguishability of the vectors $w_{5,\;l}$, $w_{6,\;l}$,
$w_{7,\;l}$, $l\geq1$, and $w_{1,\;1}$, $w_{1,\;2}$,
we have $\delta_{1}^{2}=\varepsilon_{1}^{1}=\gamma_{1}^{1}=a_{1}$. In consequence, we obtain
$$
\begin{array}{llllllll}
w_{2,\;1}: & a_{1} & \beta_{1}^{1} & \beta_{1}^{2} & \beta_{1}^{3} & a_{2},\\
\\
w_{2,\;2}: & a_{1} & \beta_{2}^{1} & a_{1} & \beta_{1}^{3} & a_{2},\\
\\
w_{3,\;1}: & a_{1} & \beta_{1}^{1} & \beta_{2}^{2} & a_{1} & a_{2},\\
\\
w_{4,\;1}: & a_{1} & a_{1} & \beta_{1}^{2} & \beta_{2}^{3} & a_{2},\\
\\
w_{5,\;1}: & \eta_{1}^{1} & \beta_{1}^{1} & \beta_{2}^{2} & a_{2} & \eta_{2}^{2}. 
\end{array}
$$
Assume that $\eta_{1}^{1}\neq a_{1}$, as if $\eta_{1}^{1}=a_{1}$ and for example $\eta_{2}^{2}\neq a_{2}$,
then we would change the order of the first and fifth coordinates,  and replace the code $\varepsilon$ by $\varepsilon':=(\tau\circ\varepsilon_5,\varepsilon_2,\varepsilon_3,\varepsilon_4, \tau\circ\varepsilon_1)$, where $\tau=(1\, 2)$. 
Consider in addition all  vectors with codes
$(l,1,2,2,2)$, $l\geq2$. They can be written in the form
$$
\begin{array}{llllllr}
w_{8,\;l-1}: & \eta_{l}^{1} & \beta_{1}^{1} & \beta_{2}^{2} & a_{2} & \rho_{2}^{l-1}, & \qquad l\geq2. 
\end{array}
$$
If $\rho_{2}^{i}=\eta_{2}^{2}$ for all $i\geq1$, then the vectors $w_{5,\;1}$ and
$w_{8,\;l}$, $l\geq1,$ determine a column. Therefore, we can assume that
$\rho_{2}^{1}\neq\eta_{2}^{2}$, as if $\rho_{2}^{1}=\eta_{2}^{2}$ and for example
$\rho_{2}^{2}\neq\eta_{2}^{2}$, then we would change the code replacing $\varepsilon_1$ by $\varepsilon'_1$, defined so that $\varepsilon'_1$ restricted to $\er\setminus (\eta^1_1+\zet)$ coincides with $\varepsilon_1$, while $\varepsilon'_1$ restricted to $\eta^1_1+\zet$ coincides with the composite $(2\,3)\circ\varepsilon_1$.
Moreover, we can also assume that
$\eta_{2}^{2}\neq a_{2}$, as if $\eta_{2}^{2}=a_{2}$
and $\rho_{2}^{1}\neq a_{2}$, then we would change the code 
similarly as above replacing $\varepsilon_1$ appropriately. 
Now, let us consider all  vectors with codes
$(1,1,2,2,1)$, $(1,1,2,2,l)$, $l\geq3$.
They have the form
$$
\begin{array}{llllllr}
w_{9,\;1}: & \sigma_{1}^{1} & \beta_{1}^{1} & \beta_{2}^{2} & a_{2} & \eta_{1}^{2},\\
\\
w_{9,\;l-1}: & \sigma_{1}^{l-1} & \beta_{1}^{1} & \beta_{2}^{2} & a_{2} & \eta_{l}^{2}, & \qquad l\geq3. 
\end{array}
$$
Since $\rho_{2}^{1}\neq\eta_{2}^{2}$,  by the distinguishability of $w_{8,\;1}$ and $w_{9,\;l}$, $l\geq1,$ we have  $\sigma_{1}^{i}=\eta_{1}^{1}$
for all $i\geq1$. Thus, the vectors $w_{5,\;1}$ and
$w_{9,\;l}$, $l\geq1,$ determine a column.
\hfill$\square$

\begin{lem}\label{le6}
If a set $T$ determining a cube tiling of $\mathbb R^{6}$ contains vectors of the form$$
\begin{array}{lllllllll}
w_{1,\;1}: & a_{1} & a_{1} & a_{1} & a_{1} & a_{1} & a_{1},\\
\\
w_{1,\;2}: & a_{2} & b_{1} & a_{1} & a_{1} & a_{1} & a_{1},\\
\\
w_{1,\;3}: & a_{1} & a_{2} & b_{1} & a_{1} & a_{1} & a_{1},\\
\\
w_{1,\;4}: & b_{1} & a_{1} & a_{2} & a_{1} & a_{1} & a_{1}, 
\end{array}
$$
then it contains vectors which determine a column.
\end{lem}

\noindent
{\bf Proof}.
Consider all  vectors with codes $(2,l,1,1,1,1)$, $l\geq2$.
As  these vectors together with $w_{1,\;1}$, $w_{1,\;2}$, $w_{1,\;3}$ and $w_{1,\;4}$  are distinguishable, they can be written in the following form 
$$
\begin{array}{lllllllr}
w_{2,\;l-1}: & a_{2} & b_{l} & a_{1} & \delta_{1}^{3l-5} & \delta_{1}^{3l-4} & \delta_{1}^{3l-3}, & \qquad l\geq2. 
\end{array}
$$
If $\delta_{1}^{i}=a_{1}$ for all $i\geq1$, then the vectors $w_{1,\;2}$, $w_{2,\;l}$, $l\geq1,$
determine a column. Therefore, let us suppose that at least one of the coordinates
$\delta_{1}^{i}$, $i\geq1$, is distinct from $a_{1}$. We can assume that it is
$\delta_{1}^{1}$, as if $\delta_{1}^{1}=a_{1}$ and for example $\delta_{1}^{5}\neq a_{1}$,
then we would change the order of the fourth and fifth coordinates, and the natural code replacing $\varepsilon_{2}$ by $\varepsilon'_{2}$ defined so that $\varepsilon'_{2}$ restricted to $\er\setminus(b_{1}+\zet)$ coincides with $\varepsilon_{2}$ while $\varepsilon'_{2}$ restricted to $b_{1}+\zet$ coincides with $(2\, 3)\circ\varepsilon_{2}$.
Now, consider all  vectors with codes
$(2,2,1,l,1,1)$, $l\geq2$. They have the form
$$
\begin{array}{lllllllr}
w_{3,\;l-1}: & a_{2} & b_{2} & a_{1} & \delta_{l}^{1} & \varepsilon_{1}^{2l-3} & \varepsilon_{1}^{2l-2}, & \qquad l\geq2. 
\end{array}
$$
If $\varepsilon_{1}^{i}=\delta_{1}^{2}$ for $i=1,3,\ldots$ and
$\varepsilon_{1}^{i}=\delta_{1}^{3}$ for $i=2,4,\ldots$, then the vectors
$w_{2,\;1}$ and $w_{3,\;l}$, $l\geq1$, determine a column. Therefore, let us assume that at least one of the above equalities
does not happen. We can assume that $\varepsilon_{1}^{1}\neq\delta_{1}^{2}$, as if
$\varepsilon_{1}^{1}=\delta_{1}^{2}$ and for example $\varepsilon_{1}^{4}\neq\delta_{1}^{3}$,
then we would change the order of the fifth and sixth coordinates, and the natural code replacing $\varepsilon_{4}$ by $\varepsilon'_{4}$ so that $\varepsilon'_{4}=(2\, 3)\circ\varepsilon_{4}$ on $\delta_{1}^{1}+\zet$ and $\varepsilon'_{4}=\varepsilon_{4}$ on the complement of $\delta_{1}^{1}+\zet$.
Moreover, we can also assume that $\delta_{1}^{2}\neq a_{1}$. (If $\delta_{1}^{2}=a_{1}$ and
$\varepsilon_{1}^{1}\neq a_{1}$, then we would change  the natural code replacing $\varepsilon_{4}$ by $\varepsilon'_{4}$ so that $\varepsilon'_{4}=(1\,2)\circ\varepsilon_{4}$ on $\delta_{1}^{1}+\zet$
and $\varepsilon'_{4}=\varepsilon_{4}$ on the complement of $\delta_{1}^{1}+\zet$.)
Now, let us take all  vectors with codes
$(2,2,1,1,l,1)$, $l\geq2$. They can be written in the form
$$
\begin{array}{lllllllr}
w_{4,\;l-1}: & a_{2} & b_{2} & a_{1} & \varphi_{1}^{2l-3} & \delta_{l}^{2} & \varphi_{1}^{2l-2}, & \qquad l\geq2. 
\end{array}
$$
By the distinguishability of the vectors $w_{3,\;1}$ and $w_{4,\;l}$, $l\geq1$, we have
$\varphi_{1}^{i}=\delta_{1}^{1}$ for $i=1,3,\ldots$, as $\varepsilon_{1}^{1}\neq\delta_{1}^{2}$. If now
$\varphi_{1}^{i}=\delta_{1}^{3}$ for $i=2,4,\ldots$, then $w_{2,\;1}$, $w_{4,\;l}$, $l\geq1$,
determine a column. Therefore, let us suppose that at least one of the coordinates
$\varphi_{1}^{i}$, $i=2,4,\ldots$, is distinct from $\delta_{1}^{3}$. We can assume that this is
$\varphi_{1}^{2}$.
Moreover, we can also assume that $\delta_{1}^{3}\neq a_{1}$. (If $\delta_{1}^{3}=a_{1}$,
and $\varphi_{1}^{2}\neq a_{1}$, then we would replace  $\varepsilon_{5}$  by $\varepsilon'_{5}$ so that $\varepsilon'_{5}=(1\,2)\circ\varepsilon_{5}$ on $\delta_{1}^{2}+\zet$ and $\varepsilon'_5=\varepsilon_{5}$ on the complement of the set $\delta_{1}^{2}+\zet$.)
Consider all  vectors with codes
$(2,2,1,1,1,l)$, $l\geq 2$. By the distinguishability, they have the form
$$
\begin{array}{lllllllr}
w_{5,\;l-1}: & a_{2} & b_{2} & a_{1} & \mu_{1}^{2l-3} & \mu_{1}^{2l-2} & \delta_{l}^{3},  & \qquad l\geq2. 
\end{array}
$$
Since $\varphi_{1}^{2}\neq\delta_{1}^{3}$, by the distinguishability of $w_{4,\;1}$ and $w_{5,\;l}$, $l\geq1$, we have
$\mu_{1}^{i}=\delta_{1}^{2}$ for $i=2,4,\ldots$. Since the vectors $w_{3,\;1}$
and $w_{5,\;l}$, $l\geq1,$ must be distinguishable, we have $\mu_{1}^{i}=\delta_{1}^{1}$
for $i=1,3,\ldots$ or $\varepsilon_{1}^{2}=\delta_{1}^{3}$.
If the first case takes place, then the vectors $w_{2,\;1}$, $w_{5,\;l}$, $l\geq1,$
determine a column. Hence, the second case must be considered. Then we obtain
$$
\begin{array}{lllllllll}
w_{2,\;1}: & a_{2} & b_{2} & a_{1} & \delta_{1}^{1} & \delta_{1}^{2} & \delta_{1}^{3},\\
\\
w_{3,\;1}: & a_{2} & b_{2} & a_{1} & \delta_{2}^{1} & \varepsilon_{1}^{1} & \delta_{1}^{3},\\
\\
w_{4,\;1}: & a_{2} & b_{2} & a_{1} & \delta_{1}^{1} & \delta_{2}^{2} & \varphi_{1}^{2},\\
\\
w_{5,\;1}: & a_{2} & b_{2} & a_{1} & \mu_{1}^{1} & \delta_{1}^{2} & \delta_{2}^{3}.
\end{array}
$$
Now, let us take into account all  vectors with codes $(2,2,1,2,l,1)$, $l\geq2$.
By their distinguishability from $w_{2,\;1}$, $w_{3,\;1}$, $w_{4,\;1}$ and $w_{5,\;1}$, they have the following form
$$
\begin{array}{lllllllr}
w_{6,\;l-1}: & \eta_{2}^{3l-5} & \eta_{2}^{3l-4} & \eta_{1}^{3l-3} & \delta_{2}^{1} & \varepsilon_{l}^{1} & \delta_{1}^{3}, & \qquad l\geq2. 
\end{array}
$$
If $\eta_{2}^{i}=a_{2}$ for $i=1,4\ldots$, $\eta_{2}^{i}=b_{2}$ for $i=2,5,\ldots$
and $\eta_{1}^{i}=a_{1}$ for $i=3,6,\ldots$, then the vectors $w_{3,\;1}$, $w_{6,\;l}$, $l\geq1$, determine a column.
Therefore, let us assume that at least one of the above equalities does not happen. Then,
by the distinguishability 
of $w_{6,\;l}$, $l\geq1$, from $w_{1,\;1}$, $w_{1,\;2}$, $w_{1,\;3}$, $w_{1,\;4}$ we obtain
$\varepsilon_{1}^{1}=a_{1}$. Similarly, taking the vectors with codes $(2,2,1,1,2,l)$ and $(2,2,1,l,1,2)$, $l\geq2$, we can show  $\varphi_{1}^{2}=a_{1}$
and $\mu_{1}^{1}=a_{1}$.
As a result, we have
$$
\begin{array}{lllllllll}
w_{3,\;1}: & a_{2} & b_{2} & a_{1} & \delta_{2}^{1} & a_{1} & \delta_{1}^{3},\\
\\
w_{4,\;1}: & a_{2} & b_{2} & a_{1} & \delta_{1}^{1} & \delta_{2}^{2} & a_{1},\\
\\
w_{5,\;1}: & a_{2} & b_{2} & a_{1} & a_{1} & \delta_{1}^{2} & \delta_{2}^{3}. 
\end{array}
$$
If we repeat the above reasoning starting with vectors with codes $(1,2,l,1,1,1)$, $l\geq2$, instead of $(2,l,1,1,1,1)$, then we can add up the following vectors
$$
\begin{array}{lllllllll}
w_{7,1}: & a_{1} & a_{2} & b_{2} & \lambda_{2}^{1} & a_{1} & \lambda_{1}^{3},\\
\\
w_{8,1}: & a_{1} & a_{2} & b_{2} & \lambda_{1}^{1} & \lambda_{2}^{2} & a_{1},\\
\\
w_{9,1}: & a_{1} & a_{2} & b_{2} & a_{1} & \lambda_{1}^{2} & \lambda_{2}^{3}.
\end{array}
$$
However, it can happen that in order to get such a system the change of the order of fourth, fifth and sixth coordinates
as well as the change of the code are to be performed.
Such changes can affect the vectors $w_{3,\;1}$, $w_{4,\;1}$, $w_{5,\;1}$ but not in a substantial way. (The scheme remains unchanged.)

Now, we repeat the whole procedure once again
beginning from the vectors with codes $(l, 1, 2, 1, 1, 1)$, $l\geq2$. Then we can add up the following vectors
$$
\begin{array}{lllllllll}
w_{10,1}: & b_{2} & a_{1} & a_{2} & \nu_{2}^{1} & a_{1} & \nu_{1}^{3},\\
\\
w_{11,1}: & b_{2} & a_{1} & a_{2} & \nu_{1}^{1} & \nu_{2}^{2} & a_{1},\\
\\
w_{12,1}: & b_{2} & a_{1} & a_{2} & a_{1} & \nu_{1}^{2} & \nu_{2}^{3}.
\end{array}
$$
As previously, the vectors $w_{3,\;1}$, $w_{4,\;1}$, $w_{5,\;1}$ and $w_{7,\;1}$, $w_{8,\;1}$, $w_{9,\;1}$
can be affected. The resulting system can be written as follows
$$
\begin{array}{lllllllll}
w_{3,\;1}: & a_{2} & b_{2} & a_{1} & A_{2}^{1} & a_{1} & A_{1}^{3},\\
\\
w_{4,\;1}: & a_{2} & b_{2} & a_{1} & A_{1}^{1} & A_{2}^{2} & a_{1},\\
\\
w_{5,\;1}: & a_{2} & b_{2} & a_{1} & a_{1} & A_{1}^{2} & A_{2}^{3},\\
\\
w_{7,\;1}: & a_{1} & a_{2} & b_{2} & B_{2}^{1} & a_{1} & B_{1}^{3},\\
\\
w_{8,\;1}: & a_{1} & a_{2} & b_{2} & B_{1}^{1} & B_{2}^{2} & a_{1},\\
\\
w_{9,\;1}: & a_{1} & a_{2} & b_{2} & a_{1} & B_{1}^{2} & B_{2}^{3},\\
\\
w_{10,\;1}: & b_{2} & a_{1} & a_{2} & \nu_{2}^{1} & a_{1} & \nu_{1}^{3},\\
\\
w_{11,\;1}: & b_{2} & a_{1} & a_{2} & \nu_{1}^{1} & \nu_{2}^{2} & a_{1},\\
\\
w_{12,\;1}: & b_{2} & a_{1} & a_{2} & a_{1} & \nu_{1}^{2} & \nu_{2}^{3},
\end{array}
$$
where capital Roman letters with lower and upper indices have a similar meaning to that of Greek letters
except that the lower index may not be equal to the code value. Precisely, they are different from all $a_i$, $i\ge 1$;  any two symbols that differ only by lower indices  represent reals whose distance is a non-zero integer; e.g. $A_{1}^{1}-A_{2}^{1}\in \zet\setminus\{0\}$.     

Now, let us consider all  vectors with codes $(2,1,1,l,1,1)$, $l\geq2$.
Since $w_{1,\;2}$ has its code equal to $(2,1,1,1,1,1)$, these vectors are as follows
$$
\begin{array}{lllllllr}
w_{13,\;l-1}: & \rho_{2}^{5l-9} & \rho_{1}^{5l-8} & \rho_{1}^{5l-7} & a_{l} & \rho_{1}^{5l-6} & \rho_{1}^{5l-5}, & \qquad l\geq2. 
\end{array}
$$
Similarly, the vectors with codes $(2,1,1,1,l,1)$, $l\geq2,$
and $(2,1,1,1,1,l)$, $l\geq 2$, can be written as follows
$$
\begin{array}{lllllllr}
w_{14,\,l-1}: & \sigma_{2}^{5l-9} & \sigma_{1}^{5l-8} & \sigma_{1}^{5l-7} & \sigma_{1}^{5l-6} & a_{l} & \sigma_{1}^{5l-5},\\
\\
w_{15,\,l-1}: & \tau_{2}^{5l-9} & \tau_{1}^{5l-8} & \tau_{1}^{5l-7} & \tau_{1}^{5l-6} & \tau_{1}^{5l-5} & a_{l}, & \qquad l\geq2. 
\end{array}
$$
If $\rho_{1}^{i}=a_{1}$ for $i=4,9,\ldots$ and $i=5,10,\ldots$, then
by the distinguishability  of $w_{13,\;l}$, $l\geq1$, from $w_{10,\;1}$, $w_{3,\;1}$ and
$w_{7,\;1}$, we have $\rho_{2}^{i}=a_{2}$
for $i=1,6,\ldots$, $\rho_{1}^{i}=b_{1}$ for $i=2,7,\ldots$ and
$\rho_{1}^{i}=a_{1}$ for $i=3,8,\ldots$. Then $w_{1,\;2}$, $w_{13,\;l}$, $l\geq1,$
determine a column. If $\sigma_{1}^{i}=a_{1}$ for $i=4,9,\ldots$ and $i=5,10,\ldots$,
or $\tau_{1}^{i}=a_{1}$ for $i=4,9,\ldots$ and $i=5,10,\ldots$, then, in the same way,
we show that the set $T$ contains vectors determining a column. Suppose  that none of these situations happen.   We can assume that $\rho_{1}^{4}\neq a_{1}$. (If for example $\rho_{1}^{10}\neq a_{1}$ and
$\rho_{1}^{4}=a_{1}$, then we would change the order of the fifth and sixth coordinates, and replace $\varepsilon_{4}$ by $\varepsilon'_{4}$ so that
$\varepsilon'_{4}=(2\,3)\circ\varepsilon_{4}$ on $a_{1}+\zet$ and $\varepsilon'_{4}=\varepsilon_{4}$ on the complement of $a_{1}+\zet$.)
Then by the distinguishability of
$w_{13,\;1}$ and $w_{14,\;l}$, $l\geq1$, we have
$\sigma_{1}^{i}=a_{1}$ for $i=4,9,\ldots$. We can assume that $\sigma_{1}^{5}\neq a_{1}$, as $\sigma_{1}^{i}\neq a_{1}$ for some $i\in\{5,10,\ldots\}$.
By the distinguishability of $w_{14,\;1}$ and $w_{15,\;l}$, $l\geq1$,
we have $\tau_{1}^{i}=a_{1}$ for $i=5,10,\ldots$. Then we can assume that
$\tau_{1}^{4}\neq a_{1}$. 
By the distinguishability of $w_{15,\;1}$ and
$w_{13,\;l}$, $l\geq1$, we have $\rho_{1}^{i}=a_{1}$ for $i=5,10,\ldots$.
As a result, we obtain
$$
\begin{array}{lllllllll}
w_{13,\;1}: & \rho_{2}^{1} & \rho_{1}^{2} & \rho_{1}^{3} & a_{2} & \rho_{1}^{4} & a_{1},\\
\\
w_{14,\;1}: & \sigma_{2}^{1} & \sigma_{1}^{2} & \sigma_{1}^{3} & a_{1} & a_{2} & \sigma_{1}^{5},\\
\\
w_{15,\;1}: & \tau_{2}^{1} & \tau_{1}^{2} & \tau_{1}^{3} & \tau_{1}^{4} & a_{1} & a_{2}.
\end{array}
$$
Now, let us take into account all  vectors with codes $(2,1,1,2,l,1)$, $l\geq2$.
By their distinguishability from $w_{13,\;1}$, $w_{14,\;1}$ and $w_{15,\;1}$, they have the following form
$$
\begin{array}{lllllllr}
w_{16,\;l-1}: & \theta_{2}^{3l-5} & \theta_{1}^{3l-4} & \theta_{1}^{3l-3} & a_{2} & \rho_{l}^{4} & a_{1}, & \qquad l\geq2. 
\end{array}
$$
By the distinguishability of $w_{13,\;1}$ and $w_{16,\;l}$, $l\geq1$, from
$w_{10,\;1}$, $w_{3,\;1}$ and
$w_{7,\;1}$,
we have $\rho_{1}^{3}=a_{1}$, $\rho_{1}^{2}=b_{1}$, $\rho_{2}^{1}=a_{2}$,
$\theta_{1}^{i}=a_{1}$ for $i=3,6,\ldots$, $\theta_{1}^{i}=b_{1}$
for $i=2,5,\ldots$ and $\theta_{2}^{i}=a_{2}$ for $i=1,4,\ldots$.
Thus, the vectors $w_{13,\;1}$ and $w_{16,\;l}$, $l\geq1,$ determine a column.
\hfill$\square$
\begin{twr}\label{ge6}
Every cube tiling of $\mathbb R^{6}$ contains a column.
\end{twr}

\noindent
{\bf Proof}.
Let $T$ be an arbitrary set which determines a cube tiling of $\mathbb R^{6}$.  By Theorems \ref{sys} and \ref{ge5}, the set $T$ contains
vectors which determine a $4$-column.
Passing to an isomorphic system if necessary, we can assume that the vectors determining our $4$-column are as follows
$$
\begin{array}{lllllllr}
w_{1,\;1}: & a_{1} & a_{1} & a_{1} & a_{1} & a_{1} & a_{1},\\
\\
w_{1,\;l}: & a_{l} & a_{1} & a_{1} & a_{1} & a_{1} & \alpha_{1}^{l-1}, & \qquad l\geq2. 
\end{array}
$$
If  $\alpha_{1}^{i}=a_{1}$ for all $i\geq1$, then the vectors $w_{1,\;k}$, $k\geq1$, determine a column.
Suppose that $\alpha_{1}^{i}\neq a_{1}$ for some $i\geq1$. We can assume that
$\alpha_{1}^{1}\neq a_{1}$. Consider all vectors with codes
$(1,\ast,\ast,\ast,\ast,l)$, $l\geq2$, where $\ast$ can take any value from
$\mathbb N$, such that their middle coordinates (second, third, fourth and fifth)
are different from $a_{l}$, $l\geq2$.
Let us pick a vector from among them whose middle coordinates differ from $a_{1}$ at as many places as possible.
Similarly as in the proof of Theorem \ref{ge4} we can assume that this vector has the code
$(1,1,1,1,1,2)$. By its distinguishability from $w_{1,1}$ and $w_{1,2}$, it can be written in the form
$$
\begin{array}{lllllll}
w_{2,\;1}: & a_{1} & \beta_{1}^{1} & \beta_{1}^{2} & \beta_{1}^{3} & \beta_{1}^{4} & a_{2}. 
\end{array}
$$
Five cases have to be considered:

\medskip
\noindent
\textit{Case 1.}
$\beta_{1}^{1}=\beta_{1}^{2}=\beta_{1}^{3}=\beta_{1}^{4}=a_{1}.$
\medskip

Then take all vectors with codes
$(1,1,1,1,1,l)$, $l\geq3$. Since each of them must be distinguishable from $w_{1,\;1}$, $w_{1,\;2}$ and $w_{2,\;1}$, they have the form
$$
\begin{array}{lllllllr}
w''''_{2,\;l-1}: & a_{1} & \beta_{1}^{4l-7} & \beta_{1}^{4l-6} & \beta_{1}^{4l-5} & \beta_{1}^{4l-4} & a_{l}, & \qquad l\geq3. 
\end{array}
$$
Since $w_{2,\;1}$ has the smallest possible number of the middle coordinates equal to $a_{1}$,
we have $\beta_{1}^{i}=a_{1}$ for all $i\geq5$. Thus, the vectors $w_{2,\;1}$, $w''''_{2,\;l}$, $l\geq2$,
determine a column.

\medskip
\noindent
\textit{Case 2.}
Exactly three of the coordinates  $\beta_{1}^{1}$, $\beta_{1}^{2}$, $\beta_{1}^{3}$, $\beta_{1}^{4}$
are equal to $a_{1}$.
\medskip

Then we can assume that  $\beta_{1}^{2}=\beta_{1}^{3}=\beta_{1}^{4}=a_{1}$,
$\beta_{1}^{1}\neq a_{1}$, as in the other case
we would change the order of the appropriate coordinates. (If $\beta_{1}^{1}=a_{1}$ and for example $\beta_{1}^{4}\neq a_{1}$,
then we exchange the second coordinate with the fifth.) Take into account all vectors with codes $(1,l,1,1,1,2)$, $l\geq2$.
By their distinguishability from $w_{1,\;1}$, $w_{1,\;2}$ and $w_{2,\;1}$, they have the following form
$$
\begin{array}{lllllllr}
w'''_{2,\;l}: & a_{1} & \beta_{l}^{1} & \gamma_{1}^{3l-5} & \gamma_{1}^{3l-4} & \gamma_{1}^{3l-3} & a_{2}, & \qquad l\geq2.
\end{array}
$$
As $w_{2,\;1}$ has the smallest possible number of the middle coordinates equal to $a_{1}$,
we have $\gamma_{1}^{i}=a_{1}$ for all $i\geq1$. Thus, the vectors $w_{2,\;1}$, $w'''_{2,\;l}$, $l\geq2$,
determine a column.

\medskip
\noindent
\textit{Case 3.}
Exactly two of the coordinates  $\beta_{1}^{1}$, $\beta_{1}^{2}$, $\beta_{1}^{3}$, $\beta_{1}^{4}$
are equal to $a_{1}$.
\medskip

Then, by the same reasoning as before,  we can assume that  $\beta_{1}^{3}=\beta_{1}^{4}=a_{1}$, 
$\beta_{1}^{1}\neq a_{1}$ and $\beta_{1}^{2}\neq a_{1}$. 
Consider all vectors with codes $(1,l,1,1,1,2)$
and  $(1,1,l,1,1,2)$, $l\geq2$. They can be written in the form
$$
\begin{array}{lllllllr}
w''_{2,\;l}: & a_{1} & \beta_{l}^{1} & \gamma_{1}^{3l-5} & \gamma_{1}^{3l-4} & \gamma_{1}^{3l-3} & a_{2},\\
\\
w''_{3,\;l-1}: & a_{1} & \delta_{1}^{3l-5} & \beta_{l}^{2} & \delta_{1}^{3l-4} & \delta_{1}^{3l-3} & a_{2}, & \qquad l\geq2. 
\end{array}
$$
Since the vectors $w''_{2,\;k}$, $k\geq2$, and $w''_{3,\;l}$, $l\geq1$,
must be distinguishable, we have $\delta_{1}^{i}=\beta_{1}^{1}$ for $i=1,4,\ldots$,
or $\gamma_{1}^{i}=\beta_{1}^{2}$ for $i=1,4,\ldots$. If the first case takes place,  
we have $\delta_{1}^{i}=a_{1}$ for
$i=2,5,\ldots$ and $i=3,6,\ldots$, as $w_{2,\;1}$ has the smallest possible number of the middle coordinates equal to $a_{1}$. Then the vectors $w_{2,\;1}$, $w''_{3,\;l}$, $l\geq1$,
determine a column. Otherwise, $\gamma_{1}^{i}=a_{1}$
for $i=2,5,\ldots$ and $i=3,6,\ldots$, and the vectors $w_{2,\;1}$, $w''_{2,\;k}$, $k\geq2$,
determine a column.

\medskip
\noindent
\textit{Case 4.}
Exactly one of the coordinates  $\beta_{1}^{1}$, $\beta_{1}^{2}$, $\beta_{1}^{3}$, $\beta_{1}^{4}$
is equal to $a_{1}$.
\medskip

Then we can suppose that  $\beta_{1}^{4}=a_{1}$, 
$\beta_{1}^{1}\neq a_{1}$, $\beta_{1}^{2}\neq a_{1}$ and $\beta_{1}^{3}\neq a_{1}$.
Consider all vectors with codes $(1,l,1,1,1,2)$, $l\geq2$.
By their distinguishability from $w_{1,\;1}$, $w_{1,\;2}$ and $w_{2,\;1}$, they can be written as follows 
$$
\begin{array}{lllllllr}
w'_{2,\;l}: & a_{1} & \beta_{l}^{1} & \gamma_{1}^{3l-5} & \gamma_{1}^{3l-4} & \gamma_{1}^{3l-3} & a_{2},  & \qquad l\geq2. 
\end{array}
$$
If $\gamma_{1}^{i}=\beta_{1}^{2}$ for $i=1,4,\ldots$ and $\gamma_{1}^{i}=\beta_{1}^{3}$
for $i=2,5,\ldots$, then we have $\gamma_{1}^{i}=a_{1}$ for $i=3,6,\ldots$, as $w_{2,\;1}$ has the smallest possible number of the middle coordinates equal to  $a_{1}$. Thus, the vectors $w_{2,\;1}$, $w'_{2,\;l}$, $l\geq2$,
determine a column. Therefore, let us suppose that at least one of the above equalities does not happen. We can assume that $\gamma_{1}^{1}\neq\beta_{1}^{2}$. (If $\gamma_{1}^{1}=\beta_{1}^{2}$ and for example $\gamma_{1}^{5}\neq\beta_{1}^{3}$, then we would change the order of the third and fourth coordinates, and the code replacing $\varepsilon_{2}$ by $\varepsilon'_{2}$ so that $\varepsilon'_{2}=(2\, 3)\circ\varepsilon_{2}$ on $\beta_{1}^{1}+\zet$ and $\varepsilon'_{2}=\varepsilon_{2}$ on the complement of $\beta^1_{1}+\zet$.)  
Take into account all vectors with codes
$(1,1,l,1,1,2)$, $l\geq2$. They have the form
$$
\begin{array}{lllllllr}
w'_{3,\;l-1}: & a_{1} & \delta_{1}^{3l-5} & \beta_{l}^{2} & \delta_{1}^{3l-4} & \delta_{1}^{3l-3} & a_{2}, & \qquad l\geq2. 
\end{array}
$$
Since  $\gamma_{1}^{1}\neq\beta_{1}^{2}$ and the vectors $w'_{2,\;2}$
and $w'_{3,\;l}$, $l\geq1$, are distinguishable, we have $\delta_{1}^{i}=\beta_{1}^{1}$
for $i=1,4,\ldots$. If now $\delta_{1}^{i}=\beta_{1}^{3}$ for
$i=2,5,\ldots$, then, since $\beta_{1}^{4}=a_{1}$, and $w_{2,\;1}$ has the smallest possible number of the middle coordinates equal to $a_{1}$ , we have $\delta_{1}^{i}=a_{1}$ for $i=3,6,\ldots$.
Thus, the vectors $w_{2,\;1}$, $w'_{3,\;l}$, $l\geq1$, determine a column. Therefore, let us assume that
$\delta_{1}^{2}\neq\beta_{1}^{3}$. (If $\delta_{1}^{2}=\beta_{1}^{3}$, and for example
$\delta_{1}^{4}\neq\beta_{1}^{3}$, then we would replace the code $\varepsilon_{3}$ by $\varepsilon'_{3}$ so that $\varepsilon'_{3}=(2\,3)\circ\varepsilon_{3}$ on $\beta_{1}^{2}+\zet$ and $\varepsilon'_{3}=\varepsilon_{3}$ on the complement of $\beta^2_{1}+\zet$.)  Take all vectors with codes
$(1,1,1,l,1,2)$, $l\geq2$. They can be written in the form
$$
\begin{array}{lllllllr}
w'_{4,\;l-1}: & a_{1} & \varepsilon_{1}^{3l-5} & \varepsilon_{1}^{3l-4} & \beta_{l}^{3} & \varepsilon_{1}^{3l-3} & a_{2}, & \qquad l\geq2. 
\end{array}
$$
Since $\delta_{1}^{2}\neq\beta_{1}^{3}$, by the distinguishability of the vectors
$w'_{3,\;1}$ and $w'_{4,\;l}$, $l\geq1$, we have
$\varepsilon_{1}^{i}=\beta_{1}^{2}$ for $i=2,5,\ldots$. Since $w'_{2,\;2}$
and $w'_{4,\;l}$, $l\geq1$, are distinguishable, we have  $\varepsilon_{1}^{i}=\beta_{1}^{1}$
for $i=1,4,\ldots$ or $\gamma_{1}^{i}=\beta_{1}^{3}$ for $i=2,5,\ldots$.
If the first possibility happens, then $\varepsilon_{1}^{i}=a_{1}$ for $i=3,6,\ldots$,
and the vectors $w_{2,\;1}$, $w'_{4,\;l}$, $l\geq1$, determine a column. Hence, $\gamma_{1}^{i}=\beta_{1}^{3}$ for $i=2,5,\ldots$.
Consequently, we obtain
$$
\begin{array}{lllllll}
w_{1,\;1}: & a_{1} & a_{1} & a_{1} & a_{1} & a_{1} & a_{1},\\
\\
w_{1,\;2}: & a_{2} & a_{1} & a_{1} & a_{1} & a_{1} & \alpha_{1}^{1},\\
\\
w_{2,\;1}: & a_{1} & \beta_{1}^{1} & \beta_{1}^{2} & \beta_{1}^{3} & a_{1} & a_{2},\\
\\
w'_{2,\;2}: & a_{1} & \beta_{2}^{1} & \gamma_{1}^{1} & \beta_{1}^{3} & \gamma_{1}^{3} & a_{2},\\
\\
w'_{3,\;1}: & a_{1} & \beta_{1}^{1} & \beta_{2}^{2} & \delta_{1}^{2} & \delta_{1}^{3} & a_{2},\\
\\
w'_{4,\;1}: & a_{1} & \varepsilon_{1}^{1} & \beta_{1}^{2} & \beta_{2}^{3} & \varepsilon_{1}^{3} & a_{2}.
\end{array}
$$
If now $\gamma_{1}^{1}\neq a_{1}$, then take all vectors with codes
$(1,2,l,1,1,2)$, $l\geq2$. They have the form
$$
\begin{array}{lllllllr}
w'_{5,\;l-1}: & a_{1} & \beta_{2}^{1} & \gamma_{l}^{1} & \beta_{1}^{3} & \varphi_{1}^{l-1} & a_{2}, & \qquad l\geq2.
\end{array}
$$
Since $\gamma_{1}^{1}\neq a_{1}$, $\beta_{1}^{3}\neq a_{1}$, and $w_{2,\;1}$ has the smallest possible number of the middle coordinates equal to $a_{1}$, we have
$\gamma_{1}^{3}=a_{1}$ and $\varphi_{1}^{i}=a_{1}$ for all $i\geq1$. Thus, the vectors 
$w'_{2,\;2}$, $w'_{5,\;l}$, $l\geq1$, determine a column. Similarly, we show that the set $T$
contains the vectors determining a column, if
$\delta_{1}^{2}\neq a_{1}$
or $\varepsilon_{1}^{1}\neq a_{1}$. Therefore,  the possibility $\gamma_{1}^{1}=\delta_{1}^{2}=\varepsilon_{1}^{1}=a_{1}$ must be considered.
Then we have
$$
\begin{array}{lllllll}
w_{2,\;1}: & a_{1} & \beta_{1}^{1} & \beta_{1}^{2} & \beta_{1}^{3} & a_{1} & a_{2},\\
\\
w'_{2,\;2}: & a_{1} & \beta_{2}^{1} & a_{1} & \beta_{1}^{3} & \gamma_{1}^{3} & a_{2},\\
\\
w'_{3,\;1}: & a_{1} & \beta_{1}^{1} & \beta_{2}^{2} & a_{1} & \delta_{1}^{3} & a_{2},\\
\\
w'_{4,\;1}: & a_{1} & a_{1} & \beta_{1}^{2} & \beta_{2}^{3} & \varepsilon_{1}^{3} & a_{2}.
\end{array}
$$
If now $\gamma_{1}^{3}=\delta_{1}^{3}=\varepsilon_{1}^{3}=a_{1}$, then we rename $\beta_{1}^{1}$, $\beta_{1}^{2}$, $\beta_{1}^{3}$ replacing them by $a_1$ and vice versa. Then we obtain
$$
\begin{array}{lllllll}
w_{2,\;1}: & a_{1} & a_1 & a_1 & a_1 & a_{1} & a_{2},\\
\\
w'_{2,\;2}: & a_{1} & a_{2} & \beta_{1}^{2} & a_{1} & a_1 & a_{2},\\
\\
w'_{3,\;1}: & a_{1} & a_{1} & a_{2} & \beta_{1}^3 & a_{1} & a_{2},\\
\\
w'_{4,\;1}: & a_{1} & \beta_{1}^{1} & a_{1} & a_{2} & a_1 & a_{2}.
\end{array}
$$
Subsequently, we change the order of the coordinates applying the cyclic permutation $(1\, 4\, 3\, 2)$, and change the code $\varepsilon_6$ by $\varepsilon'_6$ so that $\varepsilon'_6=(1\, 2)\circ \varepsilon_6$ on $a_1+\zet$ and 
$\varepsilon'_6=\varepsilon_6$ on the complement of $a_1+\zet$. 
As a result, we obtain
$$
\begin{array}{lllllll}
w_{2,\;1}: & a_{1} & a_{1} & a_{1} & a_{1} & a_{1} & a_{1},\\
\\
w'_{2,\;2}: & a_{2} & \beta_{1}^{2} & a_{1} & a_{1} & a_{1} & a_{1},\\
\\
w'_{3,\;1}: & a_{1} & a_{2} & \beta_{1}^{3} & a_{1} & a_{1} & a_{1},\\
\\
w'_{4,\;1}: & \beta_{1}^{1} & a_{1} & a_{2} & a_{1} & a_{1} & a_{1}.
\end{array}
$$
By Lemma \ref{le6}, the set $T$ contains vectors determining a column. Therefore, let us assume that $\gamma_{1}^{3}\neq a_{1}$. Consider all vectors with codes
$(1,2,1,1,l,2)$, $l\geq2$. By the distinguishability, they are as follows:
$$
\begin{array}{lllllllr}
w'_{6,\;l-1}: & a_{1} & \beta_{2}^{1} & \eta_{1}^{2l-3} & \eta_{1}^{2l-2} & \gamma_{l}^{3} & a_{2}, & \qquad l\geq2.
\end{array}
$$
If $\gamma_{1}^{3}\neq \varepsilon_{1}^{3}$, then by the distinguishability of the vectors
$w'_{6,\;l}$, $l\geq1$, and $w'_{4,\;1}$
we have $\eta_{1}^{i}=\beta_{1}^{3}$ for $i=2,4,\ldots$. We also have
$\eta_{1}^{i}=a_{1}$ for $i=1,3,\ldots$, as $w_{2}^{1}$ has the smallest possible number of the middle coordinates equal to $a_{1}$. Thus, the vectors $w'_{2,\;2}$, $w'_{6,\;l}$, $l\geq1$,
determine a column. If $\gamma_{1}^{3}=\varepsilon_{1}^{3}$, then in particular
$\varepsilon_{1}^{3}\neq a_{1}$. Take all vectors with codes
$(1,1,1,2,l,2)$, $l\geq2$. By the distinguishability, they have the form
$$
\begin{array}{lllllllr}
w'_{7,\;l-1}: & a_{1} & \mu_{1}^{2l-3} & \mu_{1}^{2l-2} & \beta_{2}^{3} & \varepsilon_{l}^{3} & a_{2}, & \qquad l\geq2.
\end{array}
$$
If $\varepsilon_{1}^{3}\neq\delta_{1}^{3}$, then by the distinguishability of the vectors
$w'_{3,\;1}$ and $w'_{7,\;l}$, $l\geq1$,
we have $\mu_{1}^{i}=\beta_{1}^{2}$ for $i=2,4,\ldots$. Then
$\mu_{1}^{i}=a_{1}$ for $i=1,3,\ldots$ and the vectors $w'_{4,\;1}$, $w'_{7,\;l}$, $l\geq1$,
determine a column. In the same way we show that the set $T$ contains vectors determining a column, if we assume that
$\delta_{1}^{3}\neq a_{1}$ or
$\varepsilon_{1}^{3}\neq a_{1}$. Now, let us assume that
$\gamma_{1}^{3}=\delta_{1}^{3}=\varepsilon_{1}^{3}\neq a_{1}$.
Let  us take into account vectors $w'_{2,\;2}$, $w'_{3,\;1}$ and $w'_{4,\;1}$. It should be clear that passing to an appropriate isomorphic system we can assume that 
$$
\begin{array}{lllllll}
w'_{2,\;2}: & a_{1} & a_{1} & a_{1} & a_{1} & a_{1} & a_{1},\\
\\
w'_{3,\;1}: & a_{1} & a_{2} & \beta_{2}^{2} & \beta_{1}^{3} & a_{1} & a_{1},\\
\\
w'_{4,\;1}: & a_{1} & \beta_{2}^{1} & \beta_{1}^{2} & a_{2} & a_{1} & a_{1}.
\end{array}
$$
Now, let us consider all vectors with codes $(1,1,l,1,1,1)$, $l\geq2$.
By their distinguishability from $w'_{2,\;2}$, $w'_{3,\;1}$ and $w'_{4,\;1}$, they can be written as follows:
$$
\begin{array}{lllllllr}
w'_{8,\;l-1}: & \nu_{1}^{3l-5} & a_{1} & a_{l} & a_{1} & \nu_{1}^{3l-4} & \nu_{1}^{3l-3}, & \qquad l\geq2. 
\end{array}
$$
If $\nu_{1}^{i}=a_{1}$ for all $i\geq1$, then the vectors $w'_{2,\;2}$, $w'_{8,\;l}$, $l\geq1$,
determine a column. Therefore, let us suppose that at least one of the coordinates $\nu_{1}^{i}$, $i\geq 1$,
is distinct from $a_{1}$. We can assume that $\nu_{1}^{1}\neq a_{1}$. (If not, then we would change the order of coordinates and the code $\varepsilon_{3}$ appropriately.)
Take into account all vectors with codes
$(l,1,2,1,1,1)$, $l\geq2$. They have the form
$$
\begin{array}{lllllllr}
w'_{9,\;l-1}: & \nu_{l}^{1} & a_{1} & a_{2} & a_{1} & \rho_{1}^{2l-3} & \rho_{1}^{2l-2}, & \qquad l\geq2. 
\end{array}
$$
If $\rho_{1}^{i}=\nu_{1}^{2}$ for $i=1,3,\ldots$ and $\rho_{1}^{i}=\nu_{1}^{3}$ for $i=2,4,\ldots$,
then the vectors $w'_{8,\;1}$, $w'_{9,\;l}$, $l\geq1$, determine a column. Therefore, let us assume that
$\rho_{1}^{1}\neq\nu_{1}^{2}$. (If not, then we would change the order of the fifth and sixth coordinates and the code $\varepsilon_{1}$ appropriately.) Then we can also assume that
$\nu_{1}^{2}\neq a_{1}$. (If $\rho_{1}^{1}\neq a_{1}$ and $\nu_{1}^{2}=a_{1}$, then we would replace the code $\varepsilon_{1}$ by $\varepsilon'_{1}$ so that $\varepsilon'_{1}=(1\,2)\circ\varepsilon_{1}$ on $\nu_{1}^{1}+\zet$ and $\varepsilon'_{1}=\varepsilon_{1}$ on the complement of $\nu_{1}^{1}+\zet$.)
Consider all vectors with codes $(1,1,2,1,l,1)$, $l\geq2$.
They are as follows
$$
\begin{array}{lllllllr}
w'_{10,\;l-1}: & \sigma_{1}^{2l-3} & a_{1} & a_{2} & a_{1} & \nu_{l}^{2} & \sigma_{1}^{2l-2}, & \qquad l\geq2.
\end{array}
$$
Since $\rho_{1}^{1}\neq\nu_{1}^{2}$ and the vectors $w'_{9,\;1}$ and
$w'_{10,\;l}$, $l\geq1$, are distinguishable, we have $\sigma_{1}^{i}=\nu_{1}^{1}$
for $i=1,3,\ldots$. If now $\sigma_{1}^{i}=\nu_{1}^{3}$ for
$i=2,4,\ldots$, then the vectors $w'_{8,\;1}$, $w'_{10,\;l}$, $l\geq1$,
determine a column. Therefore, let us assume that $\sigma_{1}^{2}\neq\nu_{1}^{3}$. (If not, then we would change the code $\varepsilon_{5}$ appropriately.)
Then we can also assume that
$\nu_{1}^{3}\neq a_{1}$. Take all vectors with codes
$(1,1,2,1,1,l)$, $l\geq2$. They can be written in the following form
$$
\begin{array}{lllllllr}
w'_{11,\;l-1}: & \tau_{1}^{2l-3} & a_{1} & a_{2} & a_{1} & \tau_{1}^{2l-2} & \nu_{l}^{3}, & \qquad l\geq2. 
\end{array}
$$
By the distinguishability of the vectors $w'_{11,\;l}$, $l\geq1$, and $w'_{10,\;1}$
we have $\tau_{1}^{i}=\nu_{1}^{2}$ for $i=2,4,\ldots$, as $\sigma_{1}^{2}\neq\nu_{1}^{3}$.
Since $w'_{9,\;1}$ and $w'_{11,\;l}$, $l\geq1$, are distinguishable,
we have  $\tau_{1}^{i}=\nu_{1}^{1}$ for $i=1,3,\ldots$ or
$\rho_{1}^{i}=\nu_{1}^{3}$ for $i=2,4,\ldots$. If the first possibility happens, then the vectors
$w'_{8,\;1}$, $w'_{11,\;l}$, $l\geq1$, determine a column. Hence, it remains to consider the second possibility .
Consequently, we obtain
$$
\begin{array}{lllllll}
w'_{8,\;1}: & \nu_{1}^{1} & a_{1} & a_{2} & a_{1} & \nu_{1}^{2} & \nu_{1}^{3},\\
\\
w'_{9,\;1}: & \nu_{2}^{1} & a_{1} & a_{2} & a_{1} & \rho_{1}^{1} & \nu_{1}^{3},\\
\\
w'_{10,\;1}: & \nu_{1}^{1} & a_{1} & a_{2} & a_{1} & \nu_{2}^{2} & \sigma_{1}^{2},\\
\\
w'_{11,\;1}: & \tau_{1}^{1} & a_{1} & a_{2} & a_{1} & \nu_{1}^{2} & \nu_{2}^{3}.
\end{array}
$$
Passing to an appropriate isomorphic system,
as  has been done before, we can assume that the latter system of vectors has the form 
$$
\begin{array}{lllllll}
w'_{8,\;1}: & a_{1} & a_{1} & a_{1} & a_{1} & a_{1} & a_{1},\\
\\
w'_{9,\;1}: & a_{2} & \rho_{1}^{1} & a_{1} & a_{1} & a_{1} & a_{1},\\
\\
w'_{10,\;1}: & a_{1} & a_{2} & \sigma_{1}^{2} & a_{1} & a_{1} & a_{1},\\
\\
w'_{11,\;1}: & \tau_{1}^{1} & a_{1} & a_{2} & a_{1} & a_{1} & a_{1}.
\end{array}
$$
By Lemma \ref{le6}, the set $T$ contains vectors determining a column.

\medskip
\noindent
\textit{Case 5.}
All of the coordinates  $\beta_{1}^{1}$, $\beta_{1}^{2}$, $\beta_{1}^{3}$, $\beta_{1}^{4}$
are distinct from $a_{1}$.
\medskip

Let us remind that
$$
\begin{array}{lllllll}
w_{1,\;1}: & a_{1} & a_{1} & a_{1} & a_{1} & a_{1} & a_{1},\\
\\
w_{1,\;2}: & a_{2} & a_{1} & a_{1} & a_{1} & a_{1} & \alpha_{1}^{1},\\
\\
w_{2,\;1}: & a_{1} & \beta_{1}^{1} & \beta_{1}^{2} & \beta_{1}^{3} & \beta_{1}^{4} & a_{2}.
\end{array}
$$
Now, let us consider all vectors with codes $(1,l,1,1,1,2)$, $(1,1,l,1,1,2)$, $(1,1,1,l,1,2)$ and $(1,1,1,1,l,2)$, $l\geq2$. By their distinguishability from
$w_{1,\;1}$, $w_{1,\;2}$ and $w_{2,\;1}$, they have the form
$$
\begin{array}{lllllllr}
w_{2,\;l}: & a_{1} & \beta_{l}^{1} & \gamma_{1}^{3l-5} & \gamma_{1}^{3l-4} & \gamma_{1}^{3l-3} & a_{2},\\
\\
w_{3,\;l-1}: & a_{1} & \delta_{1}^{3l-5} & \beta_{l}^{2} & \delta_{1}^{3l-4} & \delta_{1}^{3l-3} & a_{2},\\
\\
w_{4,\;l-1}: & a_{1} & \varepsilon_{1}^{3l-5} & \varepsilon_{1}^{3l-4} & \beta_{l}^{3} & \varepsilon_{1}^{3l-3} & a_{2},\\
\\
w_{5,\;l-1}: & a_{1} & \eta_{1}^{3l-5} & \eta_{1}^{3l-4} & \eta_{1}^{3l-3} & \beta_{l}^{4} & a_{2}, & \qquad l\geq2. 
\end{array}
$$
We prove now  that there is $\kappa\in \{\gamma, \delta,\varepsilon, \eta\}$ such that at least two of the families $\{\kappa_1^{3l-5}\colon l\ge 2\}$, $\{\kappa_1^{3l-4}\colon l\ge 2\}$, $\{\kappa_1^{3l-3}\colon l\geq 2\}$ are singletons $\{\beta_1^{i-1}\}$, where $i$ relates to the $i$-th coordinate of the vectors under consideration. The distinguishability of the vectors $w_{2,\;l}$ and $w_{3,\;l-1}$, $l\geq2$, implies 
$\gamma_{1}^{i}=\beta_{1}^{2}$ for $i=1,4,\ldots$ or $\delta_{1}^{i}=\beta_{1}^{1}$ for $i=1,4,\ldots$. We can assume that the first case takes place. (If
$\gamma_{1}^{i}\neq\beta_{1}^{2}$ for some $i\in\{1,4,\ldots\}$, then $\delta_{1}^{i}=\beta_{1}^{1}$ for $i=1,4,\ldots$ and we would change the order of the second and third coordinates.) By the distinguishability of the vectors $w_{2,\;l}$ and $w_{4,\;l-1}$, $l\geq2$, we have
$\gamma_{1}^{i}=\beta_{1}^{3}$ for $i=2,5,\ldots$ or $\varepsilon_{1}^{i}=\beta_{1}^{1}$ for $i=1,4,\ldots$.
If the first possibility happens, then each vector $w_{2,\;l}$, $l\geq 2$,  has the third coordinate equal to
$\beta_{1}^{2}$ and  the fourth equal to $\beta_{1}^{3}$.
Therefore, let us assume that $\varepsilon_{1}^{i}=\beta_{1}^{1}$ for $i=1,4,\ldots$.
Since $w_{3,\;l}$ and $w_{4,\;l}$, $l\geq1,$ are distinguishable, we have $\delta_{1}^{i}=\beta_{1}^{3}$ for $i=2,5,\ldots$
or $\varepsilon_{1}^{i}=\beta_{1}^{2}$ for $i=2,5,\ldots$. In the second case  each vector 
$w_{4,\;l}$, $l\geq1$, has its second and third coordinates equal to $\beta_{1}^{1}$ and $\beta_{1}^{2}$, respectively.
Hence $\delta_{1}^{i}=\beta_{1}^{3}$ for $i=2,5,\ldots$.
The distinguishability of $w_{3,\;l}$ and $w_{5,\;l}$, $l\geq1$, implies $\eta_{1}^{i}=\beta_{1}^{2}$ for $i=2,5,\ldots$
or $\delta_{1}^{i}=\beta_{1}^{4}$ for $i=3,6,\ldots$. If the second possibility happens, then each vector
$w_{3,\;l}$, $l\geq1$, has its fourth and fifth coordinates equal to 
$\beta_{1}^{3}$ and $\beta_{1}^{4}$, respectively. Therefore, we can assume that $\eta_{1}^{i}=\beta_{1}^{2}$ for $i=2,5,\ldots$.
Since $w_{2,\;l}$ and $w_{5,\;l-1}$, $l\geq 2$, are distinguishable, we have $\gamma_{1}^{i}=\beta_{1}^{4}$ for $i=3,6,\ldots$
or $\eta_{1}^{i}=\beta_{1}^{1}$ for $i=1,4,\ldots$.  In the first case each vector of
$w_{2,\;l}$, $l\geq2$, has its third and fifth coordinates equal to $\beta_{1}^{2}$ and $\beta_{1}^{4}$, respectively.
In the second case each vector of $w_{5,\;l}$, $l\geq 1,$ has the second coordinate equal to $\beta_{1}^{1}$ and the third equal to $\beta_{1}^{2}$.

As a result, we can assume that for the block $w_{2,\;l}$, $l\geq2$, at least two of the 
following equations hold: 

 $$
\{\gamma_1^{3l-5}\colon l\ge 2\}=\{\beta_1^{2}\},\,\, \{\gamma_1^{3l-4}\colon l\ge 2\}=\{\beta_1^{3}\},\,\, \{\gamma_1^{3l-3}\colon l\geq 2\}=\{\beta_1^{4}\}.
$$
(If it 
were for example the block $w_{4,\;l}$, $l\geq1$ instead of $w_{2,\;l}$, $l\geq2$, then we would change the order of the second and fourth coordinates.) 
We can also assume that $\gamma_{1}^{i}=\beta_{1}^{3}$
for $i=2,5,\ldots$ and $\gamma_{1}^{i}=\beta_{1}^{4}$
for $i=3,6,\ldots$. (If $\gamma_{1}^{i}=\beta_{1}^{2}$ for $i=1,4,\ldots$
and $\gamma_{1}^{i}=\beta_{1}^{3}$ for $i=2,5,\ldots$, then we would change the order of the third and fifth coordinates.
If $\gamma_{1}^{i}=\beta_{1}^{2}$ for $i=1,4,\ldots$ and $\gamma_{1}^{i}=\beta_{1}^{4}$
for $i=3,6,\ldots$, then we would change the order of the third and fourth coordinates.) Then we can also assume that
$\delta_{1}^{i}=\beta_{1}^{4}$ for $i=3,6,\ldots$. (If $\delta_{1}^{i}\neq\beta_{1}^{4}$ for some $i\in\{3,6,\ldots\}$,
but $\delta_{1}^{i}=\beta_{1}^{3}$ for $i=2,5,\ldots$, then we would change the order of the fourth and fifth coordinates.
If $\delta_{1}^{i}\neq\beta_{1}^{4}$ for some $i\in\{3,6,\ldots\}$ and $\delta_{1}^{i}\neq\beta_{1}^{3}$ for some $i\in\{2,5,\ldots\}$,
then by the distinguishability of the blocks $w_{5,\;l}$ and $w_{4,\;l}$ from $w_{3,\;l}$, $l\geq1$, we have $\eta_{1}^{i}=\beta_{1}^{2}$ for $i=2,5,\ldots$ and $\varepsilon_{1}^{i}=\beta_{1}^{2}$ for $i=2,5,\ldots$.
Moreover, since the vectors $w_{4,\;l}$ and $w_{5,\;l}$, $l\geq1$, are distinguishable, we have $\eta_{1}^{i}=\beta_{1}^{3}$ for $i=3,6,\ldots$
or $\varepsilon_{1}^{i}=\beta_{1}^{4}$ for $i=3,6,\ldots$.
If now the second possibility happens, i.e.
$\varepsilon_{1}^{i}=\beta_{1}^{4}$ for $j=3,6,\ldots$, then we would change the order of the coordinates applying the permutation
$(2\,3\, 4)$. If the first possibility takes place, i.e. $\eta_{1}^{i}=\beta_{1}^{3}$ for $i=3,6,\ldots$, then we use the permutation $(2\,3\,4\,5)$.) Consequently, we obtain
$$
\begin{array}{lllllllr}
w_{2,\;l}: & a_{1} & \beta_{l}^{1} & \gamma_{1}^{3l-5} & \beta_{1}^{3} & \beta_{1}^{4} & a_{2},\\
\\
w_{3,\;l-1}: & a_{1} & \delta_{1}^{3l-5} & \beta_{l}^{2} & \delta_{1}^{3l-4} & \beta_{1}^{4} & a_{2},\\
\\
w_{4,\;l-1}: & a_{1} & \varepsilon_{1}^{3l-5} & \varepsilon_{1}^{3l-4} & \beta_{l}^{3} & \varepsilon_{1}^{3l-3} & a_{2},\\
\\
w_{5,\;l-1}: & a_{1} & \eta_{1}^{3l-5} & \eta_{1}^{3l-4} & \eta_{1}^{3l-3} & \beta_{l}^{4} & a_{2}, & \qquad l\geq2.
\end{array}
$$
If $\gamma_{1}^{i}=\beta_{1}^{2}$ for $i=1,4,\ldots$, then the vectors
$w_{2,\;1}$, $w_{2,\;l}$, $l\geq2$, determine a column. Therefore, we can assume that at least one of the coordinates 
$\gamma_{1}^{i}$, $i=1,4,\ldots,$ is distinct from $\beta_{1}^{2}$.
We can suppose that it is $\gamma_{1}^{1}$. Then, by the distinguishability of the vectors
$w_{2,\;2}$ and $w_{3,\;l}$, $l\geq1$, we have $\delta_{1}^{i}=\beta_{1}^{1}$
for $i=1,4,\ldots$. If now $\delta_{1}^{i}=\beta_{1}^{3}$ for $i=2,5,\ldots$, then the vectors
$w_{2,\;1}$, $w_{3,\;l}$, $l\geq1$, determine a column. Therefore, we can assume that
$\delta_{1}^{2}\neq\beta_{1}^{3}$. Now, let us take into account the vectors $w_{2,\;1}$, $w_{2,\;2}$ and $w_{3,\;1}$.
Consequently, they have the following form
$$
\begin{array}{lllllll}
w_{2,\;1}: & a_{1} & \beta_{1}^{1} & \beta_{1}^{2} & \beta_{1}^{3} & \beta_{1}^{4} & a_{2},\\
\\
w_{2,\;2}: & a_{1} & \beta_{2}^{1} & \gamma_{1}^{1} & \beta_{1}^{3} & \beta_{1}^{4} & a_{2},\\
\\
w_{3,\;1}: & a_{1} & \beta_{1}^{1} & \beta_{2}^{2} & \delta_{1}^{2} & \beta_{1}^{4} & a_{2}.
\end{array}
$$
Passing to an appropriate isomorphic system, in much the same way as done before, we can assume that the latter system of vectors has the form
$$
\begin{array}{lllllll}
w_{2,\;1}: & a_{1} & a_{1} & a_{1} & a_{1} & a_{1} & a_{1},\\
\\
w_{2,\;2}: & a_{1} & a_{2} & \gamma_{1}^{1} & a_{1} & a_{1} & a_{1},\\
\\
w_{3,\;1}: & a_{1} & a_{1} & a_{2} & \delta_{1}^{2} & a_{1} & a_{1}.
\end{array}
$$
Let us consider all vectors with codes $(1,1,2,l,1,1)$, $l\geq2$.
By their distinguishability from $w_{2,\;1}$, $w_{2,\;2}$ and $w_{3,\;1}$, they are as follows
$$
\begin{array}{lllllllr}
w_{6,\;l-1}: & \rho_{1}^{3l-5} & a_{1} & a_{2} & \delta_{l}^{2} & \rho_{1}^{3l-4} & \rho_{1}^{3l-3}, & \qquad l\geq2. 
\end{array}
$$
If $\rho_{1}^{i}=a_{1}$ for all $i\geq1$, then the vectors $w_{3,\;1}$, $w_{6,\;l}$, $l\geq1$,
determine a column. Suppose that at least one of the coordinates $\rho_{1}^{i}$, $i\geq1$,
is distinct from $a_{1}$. We can assume that it is $\rho_{1}^{1}$. (If not, then we would change the order of coordinates and the code $\varepsilon_{4}$ appropriately.)
Take all vectors with codes $(l,1,2,2,1,1)$, $l\geq2$. They can be written as follows
$$
\begin{array}{lllllllr}
w_{7,\;l-1}: & \rho_{l}^{1} & a_{1} & a_{2} & \delta_{2}^{2} & \sigma_{1}^{2l-3} & \sigma_{1}^{2l-2}, & \qquad l\geq2.
\end{array}
$$
If $\sigma_{1}^{i}=\rho_{1}^{2}$ for $i=1,3,\ldots$ and $\sigma_{1}^{i}=\rho_{1}^{3}$ for $i=2,4,\ldots$,
then the vectors $w_{6,\;1}$, $w_{7,\;l}$, $l\geq1$, determine a column. Assume that $\sigma_{1}^{1}\neq\rho_{1}^{2}$.
(If not, then we would change the order of coordinates and the code $\varepsilon_{1}$ in an appropriate way.)
Then we can also assume that $\rho_{1}^{2}\neq a_{1}$. (If $\rho_{1}^{2}=a_{1}$, and
$\sigma_{1}^{1}\neq a_{1}$, then we would replace the code $\varepsilon_{1}$ by $\varepsilon'_{1}$ so that $\varepsilon'_{1}=(1\,2)\circ\varepsilon_{1}$ on $\rho_{1}^{1}+\zet$ and $\varepsilon'_{1}=\varepsilon_{1}$ on the complement of $\rho_{1}^{1}+\zet$.) 
Consider all vectors with codes
$(1,1,2,2,l,1)$, $l\geq2$. They have the form
$$
\begin{array}{lllllllr}
w_{8,\;l-1}: & \tau_{1}^{2l-3} & a_{1} & a_{2} & \delta_{2}^{2} & \rho_{l}^{2} & \tau_{1}^{2l-2}, & \qquad l\geq2.
\end{array}
$$
Since $\sigma_{1}^{1}\neq\rho_{1}^{2}$, by the distinguishability of $w_{7,\;1}$
and $w_{8,\;l}$, $l\geq1$, we have $\tau_{1}^{i}=\rho_{1}^{1}$ for $i=1,3,\ldots$.
If now $\tau_{1}^{i}=\rho_{1}^{3}$ for $i=2,4,\ldots$, then the vectors 
$w_{6,\;1}$, $w_{8,\;l}$, $l\geq1$, determine a column. Therefore, let us assume that
$\tau_{1}^{2}\neq\rho_{1}^{3}$ and $\rho_{1}^{3}\neq a_{1}$. (If not, then  we would change the code $\varepsilon_{5}$ in an appropriate way.)
Take all vectors with codes
$(1,1,2,2,1,l)$, $l\geq2$. They are as follows
$$
\begin{array}{lllllllr}
w_{9,\;l-1}: & \xi_{1}^{2l-3} & a_{1} & a_{2} & \delta_{2}^{2} & \xi_{1}^{2l-2} & \rho_{l}^{3}, & \qquad l\geq2. 
\end{array}
$$
The distinguishability of $w_{8,\;1}$ and $w_{9,\;l}$, $l\geq1$, implies
$\xi_{1}^{i}=\rho_{1}^{2}$ for $i=2,4,\ldots$,  as $\tau_{1}^{2}\neq\rho_{1}^{3}$.
Since $w_{7,\;1}$ and $w_{9,\;l}$, $l\geq1$, are distinguishable, we have 
$\xi_{1}^{i}=\rho_{1}^{1}$ for $i=1,3,\ldots$ or
$\sigma_{1}^{i}=\rho_{1}^{3}$ for $i=2,4,\ldots$. If the first possibility happens,
then the vectors $w_{6,\;1}$, $w_{9,\;l}$, $l\geq1$, determine a column.
Hence $\sigma_{1}^{i}=\rho_{1}^{3}$ for $i=2,4,\ldots$. As a result we obtain
$$
\begin{array}{lllllll}
w_{6,\;1}: & \rho_{1}^{1} & a_{1} & a_{2} & \delta_{2}^{2} & \rho_{1}^{2} & \rho_{1}^{3},\\
\\
w_{7,\;1}: & \rho_{2}^{1} & a_{1} & a_{2} & \delta_{2}^{2} & \sigma_{1}^{1} & \rho_{1}^{3},\\
\\
w_{8,\;1}: & \rho_{1}^{1} & a_{1} & a_{2} & \delta_{2}^{2} & \rho_{2}^{2} & \tau_{1}^{2},\\
\\
w_{9,\;1}: & \xi_{1}^{1} & a_{1} & a_{2} & \delta_{2}^{2} & \rho_{1}^{2} & \rho_{2}^{3}.
\end{array}
$$
Passing to an appropriate isomorphic system,
in much the same way as  done before, we can assume that the latter system of vectors has the form
$$
\begin{array}{lllllll}
w_{6,\;1}: & a_{1} & a_{1} & a_{1} & a_{1} & a_{1} & a_{1},\\
\\
w_{7,\;1}: & a_{2} & \sigma_{1}^{1} & a_{1} & a_{1} & a_{1} & a_{1},\\
\\
w_{8,\;1}: & a_{1} & a_{2} & \tau_{1}^{2} & a_{1} & a_{1} & a_{1},\\
\\
w_{9,\;1}: & \xi_{1}^{1} & a_{1} & a_{2} & a_{1} & a_{1} & a_{1}.
\end{array}
$$
By Lemma \ref{le6}, the set  $T$ contains vectors determining a column.
\hfill$\square$

\end{document}